# Convex Reconstruction of Structured Matrix Signals from Linear Measurements (I): Theoretical Results

*Yuan Tian*
Software School, Dalian University of Technology, Dalian, P.R.China, 116620
12 October, 2019

**Abstract** *We investigate the problem of reconstructing n-by-n structured matrix signal* $X=(\mathbf{x}_1,\ldots,\mathbf{x}_n)$ *via convex programming, where each column* $\mathbf{x}_j$ *is a vector of s-sparsity and all columns have the same* $l_1$*-norm. The regularizer is matrix norm* $|||X|||_1:=max_j|\mathbf{x}_j|_1$. *The contribution in this paper has two parts. The first part is about conditions for stability and robustness in signal reconstruction via solving the inf-*$|||.|||_1$ *convex programming from noise-free or noisy measurements. We establish uniform sufficient conditions which are very close to necessary conditions and non-uniform conditions are also discussed. Similar as the inf-*$l_1$ *compressive sensing theory for reconstructing vector signals, a* $|||.|||_1$*-version RIP condition is established. In addition, stronger conditions are investigated to guarantee the reconstructed signal's support stability, sign stability and approximation-error robustness(e.g., with linear convergence rate relative to any matrix norm). The second part is about conditions on number of measurements for robust reconstruction in noise. We take the convex geometric approach in random measurement setting and one of the critical ingredients in this approach is to estimate the related widths' bounds in case of Gaussian and non-Gaussian distributions. These bounds are explicitly controlled by signal's structural parameters r and s which determine matrix signal's column-wise sparsity and* $l_1$*-column-flatness respectively.*

## 1 Introduction

Compressive sensing develops effective methods to reconstruct signals accurately or approximately from accurate or noisy measurements by exploiting a priori knowledge about the signal, e.g., the signal's structural features[1-2]. So far in most investigations the signals are modeled as vectors of high ambient dimension. However, there are lots of applications in which the signals are naturally matrices or even tensors of high orders. For example, in modern radar systems, e.g., MIMO radar[3], the measurements can be naturally represented as $y_{kl} = \sum_{ij}\Phi_{kl,ij}X_{ij} + e_{kl}$ where each $y_{kl}$ is the echo sampled at specific time $k$ and specific receiver element $l$ in a linear or planar array; $\Phi_{kl,ij}$ is the coefficient of a linear processor determined by system transmitting/receiving and waveform features; $e_{kl}$ is the intensity of noise and clutter; $X_{ij}$, if nonzero, is the scattering intensity of a target detected in specific state cell $(i,j)$, e,g, a target at specific distance and radial speed, or at specific distance and direction, etc. In another class of applications related to signal sampling/ reconstruction, multivariable functions (waveforms) in a linear space spanned by given basis, e.g., $\{\psi_\mu(u)\varphi_\nu(v)\}_{\mu,\nu}$, are sampled as $s(u_i,v_j)=\sum_{\mu\nu}\psi_\mu(u_i)\varphi_\nu(v_j)\chi_{\mu,\nu}$ where $\chi_{\mu,\nu}$ are the signal's Fourier coefficients to be recovered from the samples $\{s(u_i,v_j)\}$. These are typical examples to reconstruct matrix signals and many of them can be naturally extended to even more general tensor signal models.

So far typical works on matrix signal compressive sensing include low-rank matrix recovery[1,4],

matrix completion, *Kronecker* compressive sensing[5-6], etc. Low-rank matrix recovery deals with how to reconstruct the matrix signal with sparse singular values from linear measurements using nuclear norm (sum of singular values) as the regularizer, *Kronecker* compressive sensing deals with how to reconstruct the matrix signal from matrix measurements via matrix $L_1$-norm $\sum_{ij}|X_{ij}|$ as the regularizer, dealing with the measurement operator in tensor-product form.

Matrix signals can have richer and more complicated structures than vector signals. When solving the reconstruction problem via convex programming, it is important to select the appropriate matrix norm or regularizer for specific signal structure. For example, $L_1$-norm is suitable for general sparsity, nuclear norm is suitable for singular-value-sparsity, and other regularizers are needed for more special or more fine-grained structures, e.g., column-wise sparsity, row-wise sparsity or some hybrid structure. Appropriate regularizer determines the reconstruction's performance.

**Contributions and Paper Organization** In this paper we investigate the problem of reconstructing $n$-by-$n$ matrix signal $X=(\mathbf{x}_1,\ldots,\mathbf{x}_n)$ by convex programming. Signal's structural features in concern are sparsity and flatness, i.e., each column $\mathbf{x}_j$ is a vector of $s$-sparsity and all columns have the same $l_1$-norm. Such signals naturally appear in some important applications, e.g., radar waveform space-time analysis, which will be investigated as an application in subsequent papers. The regularizer to be used is matrix norm $|||X|||_1:=max_j|\mathbf{x}_j|_1$ where $|.|_1$ is the $l_1$-norm on column vector space.

The contribution in this paper has two parts. The first part (sec.3 and sec.4) is about conditions for stability and robustness in signal reconstruction via solving the $inf$-$|||.|||_1$ convex programming from noise-free or noisy measurements. In sec. 3 we establish uniform sufficient conditions which are very close to necessary conditions and non-uniform conditions are also discussed. Similar as the $inf$-$l_1$ compressive sensing theory for reconstructing vector signals, a $|||.|||_1$-version RIP condition is investigated (theorem 3.5). In sec.4 stronger conditions are established to guarantee the reconstructed signal's support stability, sign stability and approximation-error robustness. For example, linear convergence rate relative to any matrix norm metric is established under a general condition in Theorem 4.1.

The second part in our work (sec.5 and sec.6) is to establish conditions on number of linear measurements for robust reconstruction in noise. We take the convex geometric approach[7-10] in random measurement setting and one of the critical ingredients in this approach is to estimate the related widths' bounds incase of Gaussian and non-Gaussian distributions. These bounds are explicitly controlled by signal's structural parameters $r$ and $s$ which determine matrix signal's column-wise sparsity and $l_1$-column-flatness respectively(e.g., lemma 5.1, 6.1 and 6.3).

Foundations for works in the first part is mainly general theory on convex optimization(e.g., first-order optimization conditions) in combination with the information on subdifferential $\partial|||X|||_1$, while foundations for the second part is mainly general theory on high-dimensional probability with some recent deep extensions. This paper is only focused on theoretical analysis. Algorithms, numerical investigations and applications will be the subjects in subsequent papers.

## 2 Basic Problems, Related Concepts and Fundamental Facts

**Conventions and Notations**: In this paper we only deal with vectors and matrices in real number field and only deal with square matrix signals for notation simplification, but all results are also true for rectangle matrix signals in complex field. Any vector $\mathbf{x}$ is regarded as column vector, $\mathbf{x}^T$ denotes its transpose (row vector). For a pair of vectors $\mathbf{x}$ and $\mathbf{y}$, $<\mathbf{x},\mathbf{y}>$ denotes their scalar product. For a pair of matrices X and Y, $<X,Y>$ denotes the scalar product $tr(X^TY)$. In particular, the *Frobenius* norm $<X,X>^{1/2}$ is denoted as $|X|_F$.

For a positive integer $s$, $\sum_s^{n\times n}$ denotes the set of $n$-by-$n$ matrices which column vectors are all of sparsity $s$, i.e., the number of non-zero components of each column vector is at most $s$. Let $S \equiv S_1 \cup \dots \cup S_n$ be a subset of $\{(i,j): i,j=1,\dots,n\}$ where each $S_j$ is a subset of $\{(i,j): i=1,\dots,n\}$ and its cardinality $|S_j| \le s$, $\sum_s^{n\times n}(S)$ denotes the set of $n$-by-$n$ matrices $\{M: M_{ij}=0 \text{ if } (i,j) \text{ not in } S\}$. S is called the *matrix signal's s-sparsity pattern*. Obviously $\sum_s^{n\times n}(S)$ is a linear space for given S and $\sum_s^{n\times n} = \cup_S \sum_s^{n\times n}(S)$.

A matrix $M=(\boldsymbol{m}_1,\dots,\boldsymbol{m}_n)$ is called *$l_1$-column-flat* if all its columns' $l_1$-norms $|\boldsymbol{m}_j|_1$ have the same value.

If $X_k$ is a group of random variables and $p(x)$ is some given probability distribution, then $X_k \sim^{iid} p(x)$ denotes that all these $X_k$'s are identically and independently sampled under this distribution.

## 2.1 Basic Problems

In this paper we investigate the problem of reconstructing $n$-by-$n$ matrix signal $X=(\mathbf{x}_1,\dots,\mathbf{x}_n)$ with $s$-sparse and $l_1$-flat column vectors $\mathbf{x}_1,\dots,\mathbf{x}_n$ (i.e., $|||X|||_1 = |\mathbf{x}_j|_1$ for all $j$) by solving the following convex programming problems. The regularizer is matrix norm $|||X|||_1 := max_j|\mathbf{x}_j|_1$.

*Problem* $MP^{(\alpha)}_{y,\Phi,\eta}$: $\quad inf\,|||Z|||_1 \quad \text{s.t. } Z \in R^{n\times n},\ |\mathbf{y}-\Phi(Z)|_\alpha \le \eta$ (2.1*a*)

In this setting $\mathbf{y}$ is a measurement vector in $R^m$ with some vector norm $|.|_\alpha$ defined on it, e.g., $|.|_\alpha$ being the $l_2$-norm. $\Phi: R^{n\times n} \to R^m$ is a linear operator and there is a matrix X (the real signal) satisfying $\mathbf{y}=\Phi(X)+\boldsymbol{e}$ where $|\boldsymbol{e}|_\alpha \le \eta$. In an equivalent component-wise formulation, $y_i = \langle\Phi_i, X\rangle + e_i$ where $\Phi_i \in R^{n\times n}$ for each $i=1,\dots,m$.

*Problem* $MP^{(\alpha)}_{y,A,B,\eta}$: $\quad inf\,|||Z|||_1 \quad \text{s.t. } Z \in R^{n\times n},\ |\mathbf{Y}-AZB^T|_\alpha \le \eta$ (2.1*b*)

In this setting Y is a matrix in space $R^{m\times m}$ with some matrix norm $|.|_\alpha$ defined on it, e.g., $|.|_\alpha$ being the *Frobenius*-norm. $\Phi_{A,B}: R^{n\times n} \to R^{m\times m}: Z \to AZB^T$ is a linear operator and there is a matrix signal X satisfying $Y= AXB^T+E$ and $|E|_\alpha \le \eta$. In an equivalent component-wise formulation, $y_{kl} = \langle\Phi_{kl}, X\rangle + e_{kl} = \sum_{ij}A_{ki}X_{ij}B_{lj} + e_{kl}$ where for each $1\le k, l\le m$ $\Phi_{kl}$ is a $n$-by-$n$ matrix with its $(i,j)$-entry as $A_{ki}B_{lj}$.

*Remark* 2.1 Throughout this paper we only consider the case $0 \le \eta < |\mathbf{y}|_\alpha$ for problem $MP^{(\alpha)}_{y,\Phi,\eta}$ and $0 \le \eta < |\mathbf{Y}|_\alpha$ for problem $MP^{(\alpha)}_{y,A,B,\eta}$ since otherwise the minimizers $X^*$ of these problems are trivially O.

For the above problems, we will investigate the real matrix signal X's reconstructability and approximation error where the measurement operator $\Phi$ and $\Phi_{A,B}$ (actually matrix A and B) are deterministic or at random. In some cases problem $MP^{(\alpha)}_{y,\Phi,\eta}$ and $MP^{(\alpha)}_{y,A,B,\eta}$ are equivalent each other but in other cases some specific hypothesis is only suitable to one of them, so it's appropriate to deal with them respectively.

## 2.2 Related Concepts

Some related concepts are presented in this subsection which are necessary and important to our work. For brevity all definitions are only presented in the form of vectors, however the generalization to the form of matrices is straightforward.

A cone $C$ is a subset in $R^n$ such that $tC$ is a subset of $C$ for any $t>0$. For a subset $K$ in $R^n$, its polar dual $K^* := \{y: \langle\mathbf{x},\mathbf{y}\rangle \le 0 \text{ for all } \mathbf{x} \text{ in } K\}$. $K^*$ is always a convex cone.

For a proper convex function $F(\mathbf{x})$, there are two important and related sets:

$$D(F,\mathbf{x}) := \{\mathbf{v}: F(\mathbf{x}+t\mathbf{v}) \le F(\mathbf{x}) \text{ for some } t>0\}$$

$$\partial F(\mathbf{x}) = \{\mathbf{u}: F(\mathbf{y}) \ge F(\mathbf{x}) + \langle\mathbf{y}-\mathbf{x}, \mathbf{u}\rangle \text{ for all } \mathbf{y}\}$$

and an important relation is $D(F,\mathbf{x})^* = \cup_{t>0} t\partial F(\mathbf{x})$.

Let $|.|$ be some vector norm and $|.|^*$ be its conjugate norm, i.e., $|\mathbf{u}|^* := max\{\langle\mathbf{x},\mathbf{u}\rangle: |\mathbf{x}| \le 1\}$ (e.g., $|||X|||_1^* = \sum_j |\mathbf{x}_j|_\infty$) then

$$\partial|\mathbf{x}| = \{\mathbf{u}: |\mathbf{x}| = \langle\mathbf{x},\mathbf{u}\rangle \text{ and } |\mathbf{u}|^* \le 1\} \quad (2.2)$$

Let $K$ be a cone in a vector space $L$ on which $\Phi$ is a linear operator, the minimum singular value of $\Phi$ with respect to $K$, norm $|.|_\beta$ on $L$ and norm $|.|_\alpha$ on $\Phi$'s image space is defined as

$$\lambda_{min,\alpha,\beta}(\Phi;K):=inf\{|\Phi\mathbf{u}|_\alpha: \mathbf{u} \text{ in } K \text{ and } |\mathbf{u}|_\beta=1\} \quad (2.3)$$

When both $|.|_\beta$ and $|.|_\alpha$ are $l_2$(or *Frobenius*) norms, $\lambda_{min,\alpha,\beta}(\Phi;K)$ is simply denoted as $\lambda_{min}(\Phi;K)$.

Let $K$ be a cone (not necessarily convex) in normed space $(L, |.|_\beta)$, its conic *Gaussian width* is defined as

$$w_\beta(K):=E_g[sup\{<\mathbf{g},\mathbf{u}>: \mathbf{u} \text{ in } K \text{ and } |\mathbf{u}|_\beta=1] \quad (2.4)$$

where $\mathbf{g}$ is the random vector on $L$ sampled under standard Gaussian distribution. When $|.|_\beta$ is $l_2$ or *Frobenius* norm on $L$, $w_\beta(K)$ is simply denoted as $w(K)$.

### 2.3 Fundamental Facts

Our research in the second part (sec.5 and sec.6) follows the convex geometric approach built upon a sequence of important results, which are summarized in this section as the fundamental facts. Originally these facts were presented for vector rather than matrix signals[7-10]. We re-present them for matrix signals in consistency with the form of our problems. For brevity, all facts are only presented with respect to problem $MP^{(\alpha)}_{y,\Phi,\eta}$ except for FACT 2.6.

**FACT 2.1** (1)Let $X\in R^{n\times n}$ be any matrix signal and $\mathbf{y}=\Phi(X)$, $X^*$ is the solution (minimizer) to the problem $MP^{(\alpha)}_{y,\Phi,\eta}$ where $\eta=0$, then $X^* = X$ *iff* $\ker\Phi\cap D(|||.|||_1,X)=\{O\}$.

(2) Let $X\in R^{n\times n}$ be any matrix signal and $\mathbf{y}=\Phi(X)+e$ where $|e|_\alpha\le\eta$, $X^*$ be the solution (minimizer) to the problem $MP^{(\alpha)}_{y,\Phi,\eta}$ where $\eta > 0$, $|.|_\beta$ be a norm on signal space to measure the reconstruction error, then

$$|X^*-X|_\beta \le 2\eta/\lambda_{min,\alpha,\beta}(\Phi;D(|||.|||_1,X))$$

**FACT 2.2** $K$ is a cone in $R^{n\times n}$(not necessarily convex), $\Phi$: $R^{n\times n}\to R^m$ is a linear operator with entries $\Phi_{kij} \sim^{iid} N(0,1)$, then for any $t > 0$:

$$P[\lambda_{min}(\Phi;K) \ge (m-1)^{1/2} - w(K) - t] \ge 1- exp(-t^2/2)$$

Combining these two facts, the following quite useful corollary can be obtained.

**FACT 2.3** Let X and $X^*$ be respectively the matrix signal and the solution to $MP^{(2)}_{y,\Phi,\eta}$ as specified in FACT 2.1(2), $\Phi_{kij} \sim^{iid} N(0,1)$, then for any $t$>0:

$$P[|X^*-X|_F \le 2\eta/((m-1)^{1/2} - w(D(|||.|||_1,X)) - t)_+] \ge 1 - exp(-t^2/2)$$

where $(u)_+:=max(u,0)$. In particular, when the measurement vector's dimension $m \ge w^2(D(|||.|||_1,X)) + Cw(D(|||.|||_1, X))$ where $C$ is some absolute constant, X can be reconstructed robustly with respect to the error norm $|X^*-X|_F$ with high probability by solving $MP^{(2)}_{y,\Phi,\eta}$.

**FACT 2.4** Let $F$ be any proper convex function and zero matrix is not in $\partial F(X)$, then $w_\beta^2(D(F, X)) \le E_G[inf\{|G-tV|_{\beta^*}^2: t>0, V \text{ in } \partial F(X)\}]$ where $|.|_{\beta^*}$ is the norm dual to $|.|_\beta$ and G is the random matrix with entries $G_{ij}\sim^{iid} N(0,1)$. In particular, when $|.|_\beta$ is $|.|_F$ then $w^2(D(F, X)) \le E_G[inf\{|G-tV|_F^2: t>0, V \text{ in } \partial F(X)\}]$.

This fact is useful to estimate the squared Gaussian width $w_\beta^2(D(F, X))$'s upper bound.

**FACT 2.5** Let X and $X^*$ be respectively the matrix signal and the solution to $MP^{(2)}_{y,\Phi,\eta}$ as specified in FACT 2.1(2), with the equivalent component-wise formulation $y_k=<\Phi_k,X>+e_k$, each $\Phi_k \sim^{iid} \Phi$ where $\Phi$ is a random matrix which satisfies the following conditions: (1)$E[\Phi]=0$; (2)There exists a constant $\alpha$>0 such that $\alpha \le E[<\Phi,U>]$ for all U: $|U|_F=1$; (3) There exists a constant $\sigma$>0 such that $P[|<\Phi,U>| \ge t] \le 2exp(-t^2/2\sigma^2)$. Let $\rho:= \sigma/\alpha$, then for any $t$>0:

$$P[\,|X^*-X|_F \le 2\eta/(c_1\alpha\rho^{-2}m^{1/2} - c_2\sigma w(D(|||.|||_1,X)) - \alpha t)_+\,] \ge 1 - exp(-c_3t^2)$$

where $c_1$, $c_2$, $c_3$ are absolute constants.

**FACT 2.6** $\Gamma$ is a subset in $n$-by-$n$ matrix space. Define the linear operator $\Phi_{A,B}$: $R^{n\times n} \to R^{m\times m}$: $Y=AXB^T$. In the equivalent component-wise formulation, $y_{kl}=<\Phi_{kl},X>=\sum_{ij}A_{ki}X_{ij}B_{lj}$ for each $1\le k,l\le m$, $\Phi_{kl} \sim^{iid} \Phi$ which is sampled under some given distribution. For any parameter $\xi>0$, define

$$Q_\xi(\Gamma;\Phi):=inf\{P[|<\Phi,U>|\geq\xi]: U \text{ in } \Gamma \text{ and } |U|_F=1\}$$

Furthermore, for each $1\leq k,l\leq m$, let $\varepsilon_{kl}\sim^{iid}$ Rademacher random variable $\varepsilon$ ($P[\varepsilon=\pm1]=1/2$) which are also independent of $\Phi$, and define

$$W(\Gamma;\Phi):=E_H[sup\{<H,U>: U \text{ in } \Gamma \text{ and } |U|_F=1\}] \quad \text{where } H:=m^{-1}\sum_{kl}\varepsilon_{kl}\Phi_{kl}=m^{-1}A^TEB,\ E=[\varepsilon_{kl}];$$

$$\lambda_{\min}(\Phi;\Gamma):=inf\{(\sum_{kl}|U_{kl}\Phi_{kl}|^2)^{1/2}: U \text{ in } \Gamma \text{ and } |U|_F=1\}$$

Then for any $\xi>0$ and $t>0$:

$$P[\lambda_{\min}(\Phi;\Gamma)\geq\xi m Q_{2\xi}(\Gamma;\Phi)-2W(\Gamma;\Phi)-\xi t]\geq 1-exp(-t^2/2)$$

*Remark*: In Fact 2.6 the definition of $\lambda_{\min}(\Phi;\Gamma)$ is the matrix version of that in subsection 2.2 with respect to *Frobenius* norm. The proof of FACT 2.5 and 2.6 (with respect to vector signals) can be found in [8]'s Theorem 6.3 and Proposition 5.1.

## 3 Basic Conditions on Matrix Signal Reconstruction

In this and next section we investigate sufficient and necessary conditions on the measurement operator for accurate and approximate signal reconstruction via solving problems $MP^{(\alpha)}_{y,\Phi,\eta}$ and $MP^{(\alpha)}_{y,A,B,\eta}$. For notation simplicity, we only deal with problem $MP^{(\alpha)}_{y,\Phi,\eta}$ and the formulation can be straightforwardly transformed into problem $MP^{(\alpha)}_{y,A,B,\eta}$.

We present conditions for accurate, stable and robust reconstruction respectively. As will be seen, these conditions are similar as those related to the regularizers with so-called *decomposable* subdifferentials. The vector's $l_1$-norm and matrix's nuclear norm are such examples. However, $\partial|||X|||_1$ is not even weakly-decomposable (i.e., there is no $W_0$ in $\partial|||X|||_1$ such that $<W_0, W-W_0>=0$ for all $W$ in $\partial|||X|||_1$). At first we prove a technical lemma 3.1 which describes $\partial|||X|||_1$'s structure.

**Lemma 3.1** For $n$-by-$n$ matrix $X=(\mathbf{x}_1,\ldots,\mathbf{x}_n)$ and matrix norm $|||X|||_1:=max_j|\mathbf{x}_j|_1$, the subdifferential

$$\partial|||X|||_1=\{(\lambda_1\xi_1,\ldots,\lambda_n\xi_n): \xi_j \text{ in } \partial|\mathbf{x}_j|_1 \text{ and } \lambda_j\geq0 \text{ for all } j,\ \lambda_1+\ldots+\lambda_n=1 \text{ and } \lambda_j=0 \text{ for } j: |\mathbf{x}_j|_1<max_k|\mathbf{x}_k|_1\}$$

*Proof* It's easy to verify the set $\{(\lambda_1\xi_1,\ldots,\lambda_n\xi_n): \xi_j$ in $\partial|\mathbf{x}_j|_1$ and $\lambda_j\geq0$ for all $j$, $\lambda_1+\ldots+\lambda_n=1$ and $\lambda_j=0$ for $j$: $|\mathbf{x}_j|_1<max_k|\mathbf{x}_k|_1\}$ is contained in $\partial|||X|||_1$: since for any $M\equiv(\lambda_1\xi_1,\ldots,\lambda_n\xi_n)$ in this set, we have

$$<M,X>=\sum_j\lambda_j<\xi_j,\mathbf{x}_j>=\sum_j\lambda_j|\mathbf{x}_j|_1=|||X|||_1\sum_j\lambda_j=|||X|||_1$$

and $|||.|||_1$'s conjugate norm $|||M|||_1^*=\sum_j\lambda_j|\xi_j|_\infty\leq\sum_j\lambda_j=1$, as a result M is in $\partial|||X|||_1$.

Now prove that any M in $\partial|||X|||_1$ has the form specified as a member in the above set. Let $M\equiv(\boldsymbol{\eta}_1,\ldots,\boldsymbol{\eta}_n)$, $|||Y|||_1\geq|||X|||_1+<Y-X,M>$ for all $Y\equiv(\mathbf{y}_1,\ldots,\mathbf{y}_n)$ implies:

$$max_j|\mathbf{y}_j|_1\geq max_j|\mathbf{x}_j|_1+\sum_j<\mathbf{y}_j-\mathbf{x}_j,\boldsymbol{\eta}_j> \quad (3.1)$$

Let $\boldsymbol{\eta}_j=|\boldsymbol{\eta}_j|_\infty\xi_j$ ( so $|\xi_j|_\infty=1$ if $\boldsymbol{\eta}_j\neq0$ ), then $max_j|\mathbf{y}_j|_1\geq max_j|\mathbf{x}_j|_1+\sum_j|\boldsymbol{\eta}_j|_\infty<\mathbf{y}_j-\mathbf{x}_j,\xi_j>$. For each $j$: $\boldsymbol{\eta}_j\neq0$ we can select a $i_j$ such that $|\xi_j(i_j)|=1$ and let $e^*_j$ be such a vector with component $e^*_j(i_j)=sgn\ \xi_j(i_j)$ and $e^*_j(i)=0$ for all $i\neq i_j$, then for $\mathbf{y}_j=\mathbf{x}_j+\boldsymbol{e}^*_j$, $j=1,\ldots,n$, (3.1) implies

$$1+max_j|\mathbf{x}_j|_1\geq max_j|\mathbf{y}_j|_1\geq max_j|\mathbf{x}_j|_1+\sum_j|\boldsymbol{\eta}_j|_\infty<\boldsymbol{e}^*_j,\xi_j>=max_j|\mathbf{x}_j|_1+\sum_j|\boldsymbol{\eta}_j|_\infty|\xi_j|_\infty=max_j|\mathbf{x}_j|_1+\sum_j|\boldsymbol{\eta}_j|_\infty$$

As a result $1\geq\sum_j|\boldsymbol{\eta}_j|_\infty$.

Furthermore for any given $i$, let $\mathbf{y}_j=\mathbf{x}_j$ for all $j\neq i$ and $\mathbf{y}_i$ be any vector satisfying $|\mathbf{y}_i|_1\leq|\mathbf{x}_i|_1$, then substitute these $\mathbf{y}_1,\ldots,\mathbf{y}_n$ into (3.1) we obtain

$$max_j|\mathbf{x}_j|_1\geq max_j|\mathbf{y}_j|_1\geq max_j|\mathbf{x}_j|_1+\sum_j<\mathbf{y}_j-\mathbf{x}_j,\boldsymbol{\eta}_j>=max_j|\mathbf{x}_j|_1+|\boldsymbol{\eta}_j|_\infty<\mathbf{y}_i-\mathbf{x}_i,\xi_i>$$

i.e., $<\mathbf{y}_i-\mathbf{x}_i,\xi_i>\leq0$. As a result, $<\mathbf{x}_i,\xi_i>\geq<\mathbf{y}_i,\xi_i>$ for all $\mathbf{y}_i$: $|\mathbf{y}_i|_1\leq|\mathbf{x}_i|_1$ so $<\mathbf{x}_i,\xi_i>\geq|\mathbf{x}_i|_1|\xi_i|_\infty=|\mathbf{x}_i|_1$, hence finally we get $<\mathbf{x}_i,\xi_i>=|\mathbf{x}_i|_1$. This (together with $|\xi_i|_\infty=1$) implies $\xi_i$ in $\partial|\mathbf{x}_i|_1$ if $\boldsymbol{\eta}_i\neq0$, for any $i=1,\ldots,n$.

In summary, we have so far proved that for any M in $\partial|||X|||_1$, M always has the form $(\lambda_1\xi_1,\ldots,\lambda_n\xi_n)$ where $\xi_j$ in $\partial|\mathbf{x}_j|_1$, $\lambda_j\geq0$ for all $j$ and $\lambda_1+\ldots+\lambda_n\leq1$. Since $|||X|||_1=<M,X>=\sum_j\lambda_j<\xi_j,\mathbf{x}_j>=\sum_j\lambda_j|\mathbf{x}_j|_1\leq max_j|\mathbf{x}_j|_1$ $\sum_j \lambda_j \leq |||X|||_1$, as a result $\lambda_1+\ldots+\lambda_n=1$ and $\lambda_j=0$ for $j$: $|\mathbf{x}_j|_1<max_k|\mathbf{x}_k|_1$. □

*Remark* 3.1: More explicitly, $\partial|||X|||_1 = \{(\lambda_1\xi_1,\ldots,\lambda_n\xi_n)$: $\lambda_j \geq 0$ for all $j$, $\lambda_1+\ldots+\lambda_n=1$ and $\lambda_j=0$ for $j$: $|\mathbf{x}_j|_1 < max_k|\mathbf{x}_k|_1$; $|\xi_j|_\infty \leq 1$ for all $j$ and $\xi_j(i) = sgn(X_{ij})$ for $X_{ij}\neq 0$ }.

### 3.1 Conditions on Φ For Accurate Reconstruction From Noise-free Measurements

At first we investigate the conditions for matrix signal reconstruction via solving the following problem:

*Problem* $MP_{y,\Phi,0}$: $\quad inf\ |||Z|||_1$ s.t. $Z \in R^{n\times n}$, $\mathbf{y} = \Phi(Z)$ (3.2)

**Theorem 3.1** Given positive integer $s$ and linear operator Φ, the signal $X\in \sum_s^{n\times n}$ is always the unique minimizer of problem $MP_{y,\Phi,0}$ where $\mathbf{y}=\Phi(X)$ if and only if

$$|||H_S|||_1 < |||H_{\sim S}|||_1 \quad (3.3)$$

for any $H=(\boldsymbol{h}_1,\ldots,\boldsymbol{h}_n) \in \ker\Phi\backslash\{O\}$ and any $s$-sparsity pattern S.

*Proof* To prove the necessity, let S be a $s$-sparsity pattern and $H \in \ker\Phi\backslash\{O\}$. Set $\mathbf{y} \equiv \Phi(H_S) = \Phi(-H_{\sim S})$ and $H_S \in \sum_s^{n\times n}$, $H_S$ should be the unique minimizer of $MP_{y,\Phi,0}$ with $-H_{\sim S}$ as its feasible solution, hence $|||H_S|||_1 < |||H_{\sim S}|||_1$.

Now prove the sufficiency. Let $X=(\mathbf{x}_1,\ldots,\mathbf{x}_n)$ be a matrix signal with its support $S = S_1\cup\ldots\cup S_n$ as a $s$-sparsity pattern (where $S_j = supp(\mathbf{x}_j)$) and let $\mathbf{y} = \Phi(X)$ . For any feasible solution Z ($\neq$X) of $MP_{y,\Phi,0}$, obviously there exists $H=(\boldsymbol{h}_1,\ldots,\boldsymbol{h}_n)$ in $\ker\Phi\backslash\{O\}$ such that Z=X+H. Since $\partial|||Z|||_1 \geq \partial|||X|||_1 + \langle H,M\rangle$ for any M in $\partial|||X|||_1$, we have

$\partial|||Z|||_1 - \partial|||X|||_1 \geq \sup\{\langle H,M\rangle$ : for any M in $\partial|||X|||_1$ }

$= \sup\{\langle H,M\rangle$ : M=E+V where $E=(\lambda_1 sgn(\mathbf{x}_1)\ldots,\lambda_n sgn(\mathbf{x}_n))$ and $V=(\lambda_1\xi_1,\ldots,\lambda_n\xi_n)$, $|\xi_j|_\infty \leq 1$,

$\lambda_j \geq 0$ for all $j$, $\lambda_1+\ldots+\lambda_n=1$ } ( by lemma 3.1 and notice $supp(sgn(\mathbf{x}_j)) = S_j = \sim supp(\xi_j)$ )

$\geq \sup\{ -|\langle H,E\rangle| + \langle H,V\rangle$: E and V specified as the above }

$= \sup\{ -|\sum_{j=1}^n \lambda_j\langle \boldsymbol{h}_{j|Sj}, sgn(\mathbf{x}_j)\rangle| + \sum_{j=1}^n \lambda_j\langle \boldsymbol{h}_{j|\sim Sj}, \xi_j\rangle$ : $\lambda_j$ and $\xi_j$ specified as the above }

$\geq -\sup|\sum_{j=1}^n \lambda_j\langle \boldsymbol{h}_{j|Sj}, sgn(\mathbf{x}_j)\rangle|$: $\lambda_j \geq 0$ for all $j$, $\lambda_1+\ldots+\lambda_n=1\} + \sup\{\langle H_{\sim S},V\rangle$: $|||V|||_1^* \leq 1$ }

(note that $|||V|||_1^*=\sum_j|\lambda_j\xi_j|_\infty \leq \sum_j|\xi_j|_\infty = 1$ where $|||.|||_1^*$ is $|||.|||_1$'s conjugate norm)

$= -\sup|\sum_{j=1}^n \lambda_j\langle \boldsymbol{h}_{j|Sj}, sgn(\mathbf{x}_j)\rangle|$: $\lambda_j \geq 0$ for all $j$, $\lambda_1+\ldots+\lambda_n=1\} + |||H_{\sim S}|||_1$ (3.4)

$= -max_j|\langle \boldsymbol{h}_{j|Sj}, sgn(\mathbf{x}_j)\rangle| + |||H_{\sim S}|||_1$

$\geq -max_j|\boldsymbol{h}_{j|Sj}|_1 + |||H_{\sim S}|||_1 = -|||H_S|||_1 + |||H_{\sim S}|||_1 > 0$ under the condition (3.3). As a result, X is the unique minimizer of $MP_{y,\Phi,0}$. □

*Remark* 3.2: (3.3) provides the uniform condition for signal reconstruction which applies to all unknown column-wise sparse signals. By a similar proof, we can also obtain a (non-uniform) sufficient condition for individual signal reconstruction, namely, for given operator Φ and unknown signal X with unknown support $S=S_1\cup\ldots\cup S_n$, if there exist $\lambda_1,\ldots,\lambda_n \geq 0$: $\lambda_1+\ldots+\lambda_n=1$ such that for any $H=(\boldsymbol{h}_1,\ldots,\boldsymbol{h}_n)$ in $\ker\Phi\backslash\{O\}$ there holds:

$$|\sum_{j=1}^n \lambda_j\langle \boldsymbol{h}_{j|Sj}, sgn(\mathbf{x}_j)\rangle| < |||H_{\sim S}|||_1 \quad (3.5)$$

then X will be the unique minimizer of $MP_{y,\Phi,0}$ where $\mathbf{y} = \Phi(X)$.

On the other hand, from $MP_{y,\Phi,0}$'s first-order optimization condition, i.e., for its minimizer X there exist M in $\partial|||X|||_1$ and a multiplier vector $\mathbf{u}$ such that

$$M + \Phi^T(\mathbf{u}) = O \quad (3.6)$$

then for any $H=(\boldsymbol{h}_1,\ldots,\boldsymbol{h}_n)$ in ker Φ we have

$$0 = \langle\Phi(H),\mathbf{u}\rangle = \langle H, \Phi^T(\mathbf{u})\rangle = -\langle H, M\rangle = -\langle H,E+V\rangle$$

where $E=(\lambda_1 sgn(\mathbf{x}_1)\ldots,\lambda_n sgn(\mathbf{x}_n))$ and $V=(\lambda_1\xi_1,\ldots,\lambda_n\xi_n)$, $|\xi_j|_\infty \leq 1$, so $-\langle H,E\rangle = \langle H,V\rangle$. Note that for the left hand side $|\langle H,E\rangle| = |\sum_{j=1}^n \lambda_j\langle \boldsymbol{h}_{j|Sj}, sgn(\mathbf{x}_j)\rangle|$ and for the right hand side $|\langle H,V\rangle| \leq \sum_{j=1}^n |\lambda_j\langle \boldsymbol{h}_{j|\sim Sj}, \xi_j\rangle|$

$\leq \sum_{j=1}^{n} \lambda_j|\boldsymbol{h}_{j|\sim Sj}|_1|\xi_j|_\infty \leq max_j\,|\boldsymbol{h}_{j|\sim Sj}|_1\sum_{j=1}^{n} \lambda_j\,|\xi_j|_\infty \leq max_j\,|\boldsymbol{h}_{j|\sim Sj}|_1 = |||H_{\sim S}|||_1$, we obtain a (relatively weak) non-uniform necessary condition, namely, for given operator Φ and unknown signal X with unknown support $S=S_1\cup ... \cup S_n$, if X is the minimizer of $MP_{y,\Phi,0}$ where $\mathbf{y} = \Phi(X)$ then there exist $\lambda_j \geq 0$ for $j=1,\ldots,n$: $\lambda_1+\ldots+\lambda_n=1$ such that for any $H=(\boldsymbol{h}_1,\ldots, \boldsymbol{h}_n)$ in ker Φ there holds the inequality

$$|\sum_{j=1}^{n} \lambda_j\langle \boldsymbol{h}_{j|Sj}, sgn(\mathbf{x}_j)\rangle| \leq |||H_{\sim S}|||_1 \quad (3.7)$$

### 3.2 Conditions on Φ For Stable Reconstruction From Noise-free Measurements

Now investigate the sufficient condition for reconstructing matrix signal via solving $MP_{y,\Phi,0}$ where $\mathbf{y}=\Phi(X)$ for some signal X which is unnecessarily sparse but $l_1$-column-flat. The established condition guarantees the minimizer $X^*$ to be a good approximation to the real signal X.

**Theorem 3.2** Given positive integer $s$ and linear operator Φ with the $s$-$|||.|||_1$ *Stable Null Space Property*, i.e., there exists a constant $0 < \rho < 1$ such that

$$|||H_S|||_1 \leq \rho\,|||H_{\sim S}|||_1 \quad (3.8)$$

for any H in ker Φ and sparsity pattern $S = S_1\cup ... \cup S_n$ where $|S_j| \leq s$. Let $Z=(\mathbf{z}_1,\ldots,\mathbf{z}_n)$ be any feasible solution to problem $MP_{y,\Phi,0}$ where $\mathbf{y} = \Phi(X)$ for some signal $X = (\mathbf{x}_1,\ldots,\mathbf{x}_n)$, then

$$|||Z - X|||_1 \leq (1-\rho)^{-1}(1+\rho)(2max_j\,\sigma_s(\mathbf{x}_j)_1 + max_j\,(|\mathbf{z}_j|_1 - |\mathbf{x}_j|_1)) \quad (3.9)$$

where $\sigma_s(\mathbf{v})_1 := |\mathbf{v}|_1 - (|v(i_1)|+\ldots+|v(i_s)|)$, $v(i_1),\ldots,v(i_s)$ are $\mathbf{v}$'s $s$ components with the largest absolute values.

In particular, for the minimizer $X^*$ of $MP_{y,\Phi,0}$ where the real signal X is $l_1$-column-flat, there is the reconstruction-error bound: $$|||X^* - X|||_1 \leq 2(1-\rho)^{-1}(1+\rho)\,max_j\,\sigma_s(\mathbf{x}_j)_1 \quad (3.10)$$

The proof follows almost the same logic of proving $l_1$-min reconstruction's stability for vector signals under the $l_1$ Null Space Property assumption (e.g., see sec. 4.2 in [1]). For presentation completeness we provide the simple proof here. The basic tool is an auxiliary inequality (which unfortunately does not hold for matrix norm $|||.|||_1$): given index subset Δ and any vector **x**, **z**, then[1]

$$|(\mathbf{x} - \mathbf{z})_{\sim\Delta}|_1 \leq |\,\mathbf{z}\,|_1 - |\,\mathbf{x}\,|_1 + |\,(\mathbf{x} - \mathbf{z})_\Delta\,|_1 + 2|\,\mathbf{x}_{\sim\Delta}\,|_1 \quad (3.11)$$

*Proof of Theorem* 3.2: For any feasible solution $Z=(\mathbf{z}_1,\ldots,\mathbf{z}_n)$ to problem $MP_{y,\Phi,0}$ where $\mathbf{y}=\Phi(X)$, there is $H = (\mathbf{h}_1,\ldots,\mathbf{h}_n)$ in ker Φ such that Z = H + X. Apply (3.11) to each column vector $\mathbf{z}_j$ and $\mathbf{x}_j$ we get

$$|\,\mathbf{h}_{j|\sim Sj}\,|_1 \leq |\,\mathbf{z}_j\,|_1 - |\,\mathbf{x}_j\,|_1 + |\,\mathbf{h}_{j|Sj}\,|_1 + 2|\,\mathbf{x}_{j|\sim Sj}\,|_1$$

Hence $|||H_{\sim S}|||_1 \equiv max_j\,|\mathbf{h}_{j|\sim Sj}\,|_1 \leq max_j\,(|\mathbf{z}_j|_1 - |\mathbf{x}_j|_1) + |||\mathbf{H}_S|||_1 + 2max_j\,|\mathbf{x}_{j|\sim Sj}| \leq max_j\,(|\mathbf{z}_j|_1 - |\mathbf{x}_j|_1) + \rho|||H_{\sim S}|||_1 + 2max_j\,|\,\mathbf{x}_{j|\sim Sj}\,|_1$ ( by (3.8)), namely :

$$|||H_{\sim S}|||_1 \leq (1-\rho)^{-1}(2max_j\,|\,\mathbf{x}_{j|\sim Sj}\,|_1 + max_j\,(|\mathbf{z}_j|_1 - |\mathbf{x}_j|_1))$$

As a result $|||H|||_1 = |||H_S|||_1 + |||H_{\sim S}|||_1 \leq (1+\rho)|||H_{\sim S}|||_1 \leq (1-\rho)^{-1}(1+\rho)(2max_j\,|\,\mathbf{x}_{j|\sim Sj}\,|_1 + max_j\,(|\mathbf{z}_j|_1 - |\mathbf{x}_j|_1))$ for any $s$-sparsity pattern S, which implies (3.9) since $min_S\,max_j\,|\,\mathbf{x}_{j|\sim Sj}\,|_1 = max_j\,\sigma_s(\mathbf{x}_j)_1$.

In particular, if Z is minimizer $X^*$ and X is $l_1$-column-flat then $|\mathbf{x}_j|_1=|||X|||_1$ for any $j$ so $max_j\,(|\mathbf{x}^*_j|_1 - |\mathbf{x}_j|_1) = |||X^*|||_1 - |||X|||_1 \leq 0$ for minimizer $X^*$, which implies (3.10). □

*Remark* 3.3: For any flat and sparse signal X, condition (3.8) guarantees X can be uniquely reconstructed by solving $MP_{y,\Phi,0}$ due to Theorem 3.1, while in this case the right hand side of (3.10) is zero, i.e., this theorem is consisted with the former one. In addition, (3.10) indicates that the error for the minimizer $X^*$ to approximate the flat but non-sparse signal X is controlled column-wisely by X's non-sparsity (measured by $max_j\,\sigma_s(\mathbf{x}_j)_1$).

*Remark* 3.4: Given positive integer $s$, sparsity pattern $S = S_1\cup ... \cup S_n$ where $|S_j| \leq s$ and linear operator Φ, Let $\Phi^T$ denote the adjoint operator and $\Phi_S$ denote the operator restricted on $\sum_s^{n\times n}(S)$(so $\Phi_S(X) = \Phi(X_S)$). If $\Phi_S{}^T\Phi_S$ is a bijection and there exists a constant $0 < \rho < 1$ such that the operator norm satisfies the inequality

$$N((\Phi_S{}^T\Phi_S)^{-1}\Phi_S{}^T\Phi_{\sim S}: |||.|||_1 \rightarrow |||.|||_1) \leq \rho \quad (3.12a)$$

then for any H in ker $\Phi$, from $\Phi_S(H_S) = \Phi(H_S) = \Phi(-H_{\sim S}) = \Phi_{\sim S}(-H_{\sim S})$ we can obtain $H_S = -(\Phi_S{}^T\Phi_S)^{-1}\Phi_S{}^T\Phi_{\sim S}(H_{\sim S})$ and then $|||H_{\sim S}|||_1 \leq N((\Phi_S{}^T\Phi_S)^{-1}\Phi_S{}^T\Phi_{\sim S}: |||.|||_1 \rightarrow |||.|||_1)|||H_{\sim S}|||_1 \leq \rho|||H_{\sim S}|||_1$, therefore (3.12a) is a uniformly sufficient condition stronger than $s$-$|||.|||_1$ Stable Null Space Property (3.8), similar as Tropp's exact recovery condition for vector signal's $l_1$-min reconstruction. By operator norm's duality $N(M: |.|_\alpha \rightarrow |.|_\beta) = N(M^T: |.|^*_\beta \rightarrow |.|^*_\alpha)$ an equivalent sufficient condition is:

$$N(\Phi_{\sim S}{}^T\Phi_S(\Phi_S{}^T\Phi_S)^{-1}: |||.|||_1^* \rightarrow |||.|||_1^*) \leq \rho \quad (3.12b)$$

The condition (3.12) can be enhanced to provide more powerful results for signal reconstruction (discussed in next section). Now we conclude this subsection with a simple condition for problem $MP_{y,\Phi,0}$ and $MP_{y,\Phi,\eta}$ to guarantee their minimizers' $l_1$-column-flatness.

**Theorem 3.3**(Condition for Minimizer's $l_1$-Column-Flatness**)** Given positive integer $s$, sparsity pattern $S = S_1 \cup \ldots \cup S_n$ where $|S_j| \leq s$ and linear operator $\Phi$, let $\Phi^T$ denote the adjoint operator and $\Phi_S$ denote the operator restricted on $\sum_S^{n\times n}(S)$ (so $\Phi_S(X) = \Phi(X_S)$). If $\Phi_S{}^T(\mathbf{z})$ doesn't have any zero-column for any $\mathbf{z} \neq \mathbf{0}$, then any minimizer $X^*$ of $MP_{y,\Phi,0}$ or $MP_{y,\Phi,\eta}$ with supp($X^*$) contained in S is $l_1$-column-flat.

*Proof* Consider the problem $MP_{y,\Phi,\eta}$ : $inf\ |||Z|||_1$ s.t. $Z \in R^{n\times n}$, $|\mathbf{y}-\Phi_S(Z)|_2 \leq \eta$ at first where $\eta > 0$. For any minimizer $X^*$ of this problem with both its objective $|||.|||_1$ and constraint function $|\mathbf{y}-\Phi_S(.)|_2$ convex, according to the general convex optimization theory, there exist a positive multiplier $\gamma^* > 0$ and $M^*$ in $\partial|||X^*|||_1$ such that

$$M^* + \gamma^*\Phi_S{}^T(\Phi_S(X^*) - \mathbf{y}) = O \text{ and } |\mathbf{y} - \Phi_S(X^*)|_2 = \eta \quad (3.13)$$

then $M^* = \gamma^*\Phi_S{}^T(\mathbf{y} - \Phi_S(X^*))$ can not have any zero column since $\mathbf{y}-\Phi_S(X^*) \neq \mathbf{0}$, which implies $|\mathbf{x}^*_j|_1 = max_k|\mathbf{x}_k{}^*|_1$ for every $j$ according to lemma 3.1.

Now consider the problem $MP_{y,\Phi,0}$ : $inf\ |||Z|||_1$ s.t. $Z \in R^{n\times n}$, $\mathbf{y} = \Phi_S(Z)$. For its minimizer $X^*$ there is a multiplier vector $\mathbf{u}$ such that $M^* + \Phi_S{}^T(\mathbf{u}) = O$. If $\mathbf{u} \neq \mathbf{0}$ then $M^*$ doesn't have any zero column which implies $|\mathbf{x}_j|_1 = max_k|\mathbf{x}_k|_1$ for every $j$ according to lemma 3.1. On the other hand, $\mathbf{u} = \mathbf{0}$ implies $M^* = O$ which cannot happen according to lemma 3.1 unless $X^* = O$. □

Note: For $MP_{y,\Phi,\eta}$ where $\Phi(Z)_k = tr(\Phi_k{}^T Z)$, the adjoint operator $\Phi^T(\mathbf{z}) = \sum_{k=1}^m z_k\Phi_k: R^m \rightarrow R^{n\times n}$. For problem $MP_{Y,A,B,\eta}$ where $\Phi_{A,B}(Z) = AZB^T$, the adjoint operator $\Phi^T(\mathbf{Y}) = A^TYB$: $R^{m\times m} \rightarrow R^{n\times n}$. In addition, supp($\Phi_S{}^T(\mathbf{z})$) is always a subset in S.

*Remark* 3.5: For the unconstrained convex optimization problem

$$X^* = Arg\ inf\ |||Z|||_1 + (1/2)\gamma\ |\mathbf{y} - \Phi_S(Z)|_2{}^2$$

The sufficient condition for minimizer $X^*$'s $l_1$-column-flatness is the same: $\Phi_S{}^T(\mathbf{z})$ doesn't have any zero-column for $\mathbf{z} \neq \mathbf{0}$. This fact will be used in sec.4.1.

In fact, the first-order optimization condition guarantees there is $M^*$ in $\partial|||X^*|||_1$ such that $M^* + \gamma\Phi_S{}^T(\Phi_S(X^*) - \mathbf{y}) = O$, i.e., $M^* = \gamma\Phi_S{}^T(\mathbf{y} - \Phi_S(X^*))$ in $\partial|||X^*|||_1$. Under the above condition, $M^*$ has no $\mathbf{0}$-column unless $\mathbf{y} - \Phi_S(X^*) = \mathbf{0}$. However, $\mathbf{y} - \Phi_S(X^*) = \mathbf{0}$ implies $M^* = O$ which cannot happen in $\partial|||X^*|||_1$ unless $X^* = O$. As a result, $|\mathbf{x}^*_j|_1 = max_k|\mathbf{x}^*_k|_1 = |||X^*|||_1$ for every $j$.

### 3.3 Conditions on Φ For Robust Reconstruction From Noisy Measurements

Now consider matrix signal reconstruction from noisy measurements by solving the convex optimization problem $MP_{y,\Phi,\eta}$ : $inf\ |||Z|||_1$ s.t. $Z \in R^{n\times n}$, $|\mathbf{y}-\Phi(Z)|_2 \leq \eta$ where $\eta > 0$.

**Theorem 3.4** Given positive integer $s$ and linear operator $\Phi$ with the $s$-$|||.|||_1$ *Robust Null Space Property*,

i.e., there exist constant $0 < \rho < 1$ and $\beta > 0$ such that

$$|||H_S|||_1 \le \rho\, |||H_{\sim S}|||_1 + \beta|\Phi(H)|_2 \quad (3.14)$$

for any $n$-by-$n$ matrix H and sparsity pattern $S = S_1 \cup \ldots \cup S_n$ where $|S_j| \le s$. Let $Z=(\mathbf{z}_1,\ldots,\mathbf{z}_n)$ be any feasible solution to the problem $MP_{y,\Phi,\eta}$ where $\mathbf{y} = \Phi(X) + e$ for some signal $X = (\mathbf{x}_1,\ldots,\mathbf{x}_n)$ and $|e|_2 \le \eta$, then

$$|||Z - X|||_1 \le (1-\rho)^{-1}(1+\rho)(2max_j\, \sigma_s(\mathbf{x}_j)_1 + 2\beta\eta + max_j\, (|\mathbf{z}_j|_1 - |\mathbf{x}_j|_1)) \quad (3.15)$$

In particular, for the minimizer $X^*$ of $MP_{y,\Phi,\eta}$ where the real signal X is $l_1$-column-flat, there is the error-control inequality:

$$|||X^* - X|||_1 \le 2(1-\rho)^{-1}(1+\rho)(\, max_j\, \sigma_s(\mathbf{x}_j)_1 + \beta\eta\,) \quad (3.16)$$

*Proof* For any feasible solution $Z=(\mathbf{z}_1,\ldots,\mathbf{z}_n)$ to problem $MP_{y,\Phi,\eta}$ where $\mathbf{y}=\Phi(X)+e$, Let $Z - X = H = (\mathbf{h}_1,\ldots,\mathbf{h}_n)$, apply (3.11) to each column vector $\mathbf{z}_j$ and $\mathbf{x}_j$ we get $|\,\mathbf{h}_{j|\sim Sj}\,|_1 \le |\,\mathbf{z}_j\,|_1 - |\,\mathbf{x}_j\,|_1 + |\,\mathbf{h}_{j|Sj}\,|_1 + 2|\,\mathbf{x}_{j|\sim Sj}\,|_1$ Hence $|||H_{\sim S}|||_1 \equiv max_j\, |\mathbf{h}_{j|\sim Sj}\,|_1 \le max_j\, (|\mathbf{z}_j|_1 - |\mathbf{x}_j|_1) + |||\mathbf{H}_S|||_1 + 2max_j\, |\mathbf{x}_{j|\sim Sj}| \le max_j\, (|\mathbf{z}_j|_1 - |\mathbf{x}_j|_1) + \rho|||H_{\sim S}|||_1 + 2max_j\, |\,\mathbf{x}_{j|\sim Sj}\,|_1 + \beta|\Phi(H)|_2$ ( by (3.14)), namely :

$$|||H_{\sim S}|||_1 \le (1-\rho)^{-1}(2max_j\, |\,\mathbf{x}_{j|\sim Sj}\,|_1 + max_j\, (|\mathbf{z}_j|_1 - |\mathbf{x}_j|_1) + \beta|\Phi(H)|_2)$$

As a result $|||H|||_1 = |||H_S|||_1 + |||H_{\sim S}|||_1 \le (1+\rho)|||H_{\sim S}|||_1 + \beta|\Phi(H)|_2 \le (1-\rho)^{-1}(1+\rho)(2max_j\, |\,\mathbf{x}_{j|\sim Sj}\,|_1 + max_j\, (|\mathbf{z}_j|_1 - |\mathbf{x}_j|_1)) + 2(1-\rho)^{-1}\beta|\Phi(X)|_2\,)$ for any $s$-sparsity pattern S, which implies (3.15) since $min_S\, max_j\, |\,\mathbf{x}_{j|\sim Sj}\,|_1 = max_j\, \sigma_s(\mathbf{x}_j)_1$.

In particular, if Z is a minimizer $X^*$ and X is $l_1$-column-flat then $|\mathbf{x}_j|_1 = |||X|||_1$ for any $j$ so $max_j\, (|\mathbf{x}^*_j|_1 - |\mathbf{x}_j|_1) = |||X^*|||_1 - |||X|||_1 \le 0$ for minimizer $X^*$, which implies (3.16). □

*Remark* 3.6: (3.16) indicates that the error for the minimizer $X^*$ to approximate the flat but non-sparse signal X is up-bounded column-wisely by X's non-sparsity (measured by $max_j\, \sigma_s(\mathbf{x}_j)_1$) and the noise strength $\eta$. If the signal X is both $s$-sparse and $l_1$-column-flat, the column-wise approximation error $|||X^* - X|||_1 \le 2(1-\rho)^{-1}(1+\rho)\beta\eta = O(\eta)$, i.e., with linear convergence rate.

*Remark* 3.7: For any minimizer $X^*$ of problem $MP_{y,\Phi,\eta}$, (3.13) is the first-order optimization condition with positive multiplier $\gamma^* > 0$ and $M^*$ in $\partial|||X^*|||_1$. Then for any $H = (\mathbf{h}_1,\ldots,\mathbf{h}_n)$ we have $\langle M^*, H\rangle = \gamma^*\langle \mathbf{y} - \Phi(X^*), \Phi(H)\rangle$ where by lemma 3.1 $M^*=E+V$, $E=(\lambda_1 sgn(\mathbf{x}_1^*)\ldots, \lambda_n sgn(\mathbf{x}_n^*))$, $V=(\lambda_1\xi_1,\ldots,\lambda_n\xi_n)$, $|\xi_j|_\infty \le 1$ and note that $supp(E) = supp(X^*) = S = \sim supp(V)$, so $\langle E,H_S\rangle = \langle E,H\rangle = -\langle V,H\rangle + \gamma^*\langle \mathbf{y} - \Phi(X^*), \Phi(H)\rangle = -\langle V,H_{\sim S}\rangle + \gamma^*\langle \mathbf{y} - \Phi(X^*), \Phi(H)\rangle$. Since for the left hand side

$$|\langle E,H_S\rangle| = |\textstyle\sum_{j=1}^{n} \lambda_j\langle \boldsymbol{h}_{j|Sj}, sgn(\mathbf{x}_j^*)\rangle\,|$$

and for the right hand side

$$|\langle V,H_{\sim S}\rangle| + |\,\gamma^*\langle \mathbf{y} - \Phi(X^*), \Phi(H)\rangle| \le \textstyle\sum_{j=1}^{n} |\lambda_j\langle \boldsymbol{h}_{j|\sim Sj}, \xi_j\rangle| + \gamma^*|\mathbf{y} - \Phi(X^*)|_2|\Phi(H)|_2$$
$$\le \textstyle\sum_{j=1}^{n} \lambda_j|\boldsymbol{h}_{j|\sim Sj}|_1|\xi_j|_\infty + \gamma^*\eta|\Phi(H)|_2 \le |||H_{\sim S}|||_1 + \gamma^*\eta|\Phi(H)|_2$$

we obtain a non-uniform necessary condition, namely, for given operator $\Phi$ and unknown signal X with unknown support $S=S_1\cup \ldots \cup S_n$ and $\mathbf{y} = \Phi(X)+\boldsymbol{e}$, if $X^*$ is the minimizer of $MP_{y,\Phi,\eta}$ with the correct support S and correct non-zero component signs as the real signal X, then there exist constants $\beta(=\gamma^*\eta) > 0$ and $\lambda_j \ge 0$ for all $j=1,\ldots,n$: $\lambda_1+\ldots+\lambda_n=1$ such that for any $H=(\boldsymbol{h}_1,\ldots, \boldsymbol{h}_n)$ there holds the inequality

$$|\textstyle\sum_{j=1}^{n} \lambda_j\langle \boldsymbol{h}_{j|Sj}, sgn(\mathbf{x}_j)\rangle\,| \le |||H_{\sim S}|||_1 + \beta|\Phi(H)|_2, \quad (3.17)$$

### 3.4 M-Restricted Isometry Property

It's well known that *RIP* with appropriate parameters provide powerful sufficient conditions to guarantee $l_1$-*min* reconstruction for sparse vector signals. With our regularizer $|||X|||_1$ we propose a similar but slightly stronger condition to guarantee reconstruction robustness by solving the convex programming $MP_{y,\Phi,\eta}$.

**Theorem 3.5** Given positive integer $s$ and linear operator $\Phi$: $R^{n\times n} \to R^m$. Suppose there exist positive constants $0 < \delta_s < 1$ and $\Delta_s > 0$ such that:

(1) $(1-\delta_s)|Z|_F^2 \le |\Phi(Z)|_2^2 \le (1+\delta_s)|Z|_F^2$ for any $Z \in \sum_s^{n\times n}$; (3.18)

(2) $|<\Phi(Z), \Phi(W)>| \le (\Delta_s/n)\sum_{j=1}^{n}|z_j|_2\,|w_j|_2$ (3.19)

for any $Z=(\mathbf{z}_1,\ldots,\mathbf{z}_n) \in \sum_s^{n\times n}$ *and* $W=(\mathbf{w}_1,\ldots,\mathbf{w}_n) \in \sum_s^{n\times n}$ with $supp(Z)\cap supp(W) = \emptyset$. Under these two conditions, there are constants $\rho$ and $\beta$ such that

$$|||H_S|||_1 \le \rho\,|||H_{\sim S}|||_1 + \beta|\Phi(H)|_2 \quad (3.20)$$

for any $n$-by-$n$ matrix H and any $s$-sparsity pattern S, where the constants can be selected as

$$\rho \le \Delta_s/(1-\delta_s-\Delta_s/4), \quad \beta \le s^{1/2}(1+\delta_s)^{1/2}/(1-\delta_s-\Delta_s/4) \quad (3.21)$$

In particular, $\delta_s + 5\Delta_s/4 < 1$ implies the robust null space condition: $\rho < 1$.

*Note*: Condition (1) is the standard RIP which implies $|<\Phi(Z),\Phi(W)>| \le \delta_{2s}|Z|_F|W|_F$ for $s$-sparse matrices Z and W with separated supports, slightly weaker than condition (2).

*Proof* Let $H = (\mathbf{h}_1,\ldots,\mathbf{h}_n)$ be any $n$-by-$n$ matrix. For each $j$ suppose $|h_j(i_1)| \ge |h_j(i_2)| \ge \ldots \ge |h_j(i_n)|$, let $S_0(j) = \{(i_1, j),\ldots,(i_s, j)\}$, i.e., the set of indices of $s$ components in column $\mathbf{h}_j$ with the largest absolute values, $S_1(j) = \{(i_{1+s}, j),\ldots,(i_{2s}, j)\}$ be the set of indices of $s$ components in $\mathbf{h}_j$ with the secondary largest absolute values, etc., and for any $k$=0,1,2,… let $S_k = \cup_{j=1}^{n} S_k(j)$, obviously $H = \sum_{k\ge 0} H_{Sk}$. At first we note that (3.20) holds for S as long as it holds for $S_0$, so we try to prove this in the following. Start from condition (1):

$(1-\delta_s)|H_{S0}|_F^2 \le |\Phi(H_{S0})|_2^2 = <\Phi(H_{S0}), \Phi(H) - \sum_{k\ge 1}\Phi(H_{Sk})>$

$= <\Phi(H_{S0}), \Phi(H)> - \sum_{k\ge 1}<\Phi(H_{S0}), \Phi(H_{Sk})>$

$\le |\Phi(H_{S0})|_2|\Phi(H)|_2 + (\Delta_s/n)\sum_{n\ge j\ge 1}\sum_{k\ge 1}|\mathbf{h}_{j|S0(j)}|_2|\mathbf{h}_{j|Sk(j)}|_2$ (by condition (2))

$\le (1+\delta_s)^{1/2}|H_{S0}|_F|\Phi(H)|_2 + (\Delta_s/n)\,|H_{S0}|_F\sum_{n\ge j\ge 1}\sum_{k\ge 1}|\mathbf{h}_{j|Sk(j)}|_2$ (by condition (1) and $|\mathbf{h}_{j|S0(j)}|_2 \le |H_{S0}|_F$)

$\le (1+\delta_s)^{1/2}|H_{S0}|_F|\Phi(H)|_2 + (\Delta_s/n)\,|H_{S0}|_F\sum_{n\ge j\ge 1}(s^{-1/2}|\mathbf{h}_{j|\sim S0(j)}|_1 + (1/4)|\mathbf{h}_{j|S0(j)}|_2)$

( by the inequality $(\sum_{s\ge k\ge 1} a_k^2)^{1/2} \le s^{-1/2}\sum_{s\ge k\ge 1} a_k + (s^{1/2}/4)(a_1 - a_s)$ for $a_1 \ge a_2 \ge \ldots \ge a_s \ge 0$

and the fact $min_{s\ge i\ge 1}|\boldsymbol{h}_{j|Sk(j)}(i)| \ge max_{s\ge i\ge 1}|\boldsymbol{h}_{j|Sk+1(j)}(i)|$ for any $j$ )

$\le |H_{S0}|_F\,((1+\delta_s)^{1/2}|\Phi(H)|_2 + (s^{-1/2}\Delta_s/n)\sum_{n\ge j\ge 1}|\mathbf{h}_{j|\sim S0(j)}|_1 + (\Delta_s/4n)\sum_{n\ge j\ge 1}|\mathbf{h}_{j|S0(j)}|_2)$

$\le |H_{S0}|_F\,((1+\delta_s)^{1/2}|\Phi(H)|_2 + s^{-1/2}\Delta_s\,max_j|\mathbf{h}_{j|\sim S0(j)}|_1 + (\Delta_s/4n)n^{1/2}(\sum_{n\ge j\ge 1}|\mathbf{h}_{j|S0(j)}|_2^2)^{1/2}$

$= |H_{S0}|_F\,((1+\delta_s)^{1/2}|\Phi(H)|_2 + s^{-1/2}\Delta_s|||H_{\sim S0}|||_1 + (\Delta_s/4n^{1/2})|H_{S0}|_F)$

Cancel $|H_{S0}|_F$ on both sides we get $(1-\delta_s)|H_{S0}|_F \le (1+\delta_s)^{1/2}|\Phi(H)|_2 + s^{-1/2}\Delta_s|||H_{\sim S0}|||_1 + (\Delta_s/4n^{1/2})|H_{S0}|_F$ hence

$$|H_{S0}|_F \le (1-\delta_s-\Delta_s/4n^{1/2})^{-1}((1+\delta_s)^{1/2}|\Phi(H)|_2 + s^{-1/2}\Delta_s|||H_{\sim S0}|||_1)$$

Note that $|||H_{S0}|||_1 = max_j|\mathbf{h}_{j|S0(j)}|_1 \le s^{1/2}max_j|\mathbf{h}_{j|S0(j)}|_2 \le s^{1/2}|H_{S0}|_F$ and combine this with the above inequality, we obtain (3.20) and (3.21) for $S_0$, which implies they hold for any S. □

## 4 More Properties on Reconstruction from Noisy Measurements

In this section we establish stronger conditions on the measurement operator $\boldsymbol{\Phi}$ for some stronger results on sparse and flat matrix signal reconstruction from noisy measurements, e.g., conditions to guarantee uniqueness, support and sign stability as well as value-error robustness.

As in last sections, for notation simplification this section only deals with problem $MP_{y,\Phi,\eta}$ but all conclusions can be straightforwardly transformed into the formulation for problem $MP_{Y,A,B,\eta}$. At first we note the basic fact that $X^* = arginf\,|||Z|||_1$ s.t. $Z \in R^{n\times n}$, $|\mathbf{y} - \Phi(Z)|_2 \le \eta$ if and only if their exists a multiplier $\gamma^* > 0$ (dependent on $X^*$ in general) such that $X^*$ is a minimizer of the unconstrained convex programming $inf\,|||Z|||_1 + (1/2)\gamma^*\,|\mathbf{y} - \Phi(Z)|_2^2$. In sec. 4.1 we investigate some critical properties for the minimizer of unconstrained optimization (4.1), then on basis of these results we establish conditions for robustness, support and sign stability in signal reconstruction via solving $MP_{y,\Phi,\eta}$ in sec. 4.2.

In the following for given positive integer $s$, sparsity pattern $S = S_1 \cup \dots \cup S_n$ where $|S_j| \le s$ for all $j$ and the linear operator $\Phi$: $R^{n\times n} \to R^m$, when $\Phi_S^T\Phi_S$ is a bijection for $\sum_s^{n\times n}(S) \to \sum_s^{n\times n}(S)$ we denote the pseudo- inverse $(\Phi_S^T\Phi_S)^{-1}\Phi_S^T$: $R^m \to \sum_s^{n\times n}(S)$ as $\Phi_S^{*-1}$.

## 4.1 Conditions on Minimizer Uniqueness and Robustness for $MLP_{y,\Phi}(\gamma)$

Consider the convex programming (4.1) with given parameter $\gamma > 0$ (value of $\gamma$ is independently set):

$$\textit{Problem } MLP_{y,\Phi}(\gamma) \qquad \inf |||Z|||_1 + (1/2)\gamma\,|\mathbf{y} - \Phi(Z)|_2^2 \tag{4.1}$$

Lemma 4.1 indicates basic properties of its sparse minimizer under some sparsity-related conditions.

**Lemma 4.1** Given **y**, positive integer $s$ and sparsity pattern $S = S_1 \cup \dots \cup S_n$ where $|S_j| \le s$ for all $j$, suppose the linear measurement operator $\Phi$: $R^{n\times n} \to R^m$ satisfies:

(1) $\Phi_S^T(\mathbf{z})$ does not have any **0**-column for $\mathbf{z} \ne \mathbf{0}$;

(2) $\Phi_S^T\Phi_S$ is a bijection;

(3) $\gamma \sup\{<\Phi_{\sim S}^T(\Phi_S\Phi_S^{*-1}(\mathbf{y}) - \mathbf{y}), H>: |||H|||_1=1\} + \sup\{<\Phi_{\sim S}^T(\Phi_S^{*-1})^T M, H>: |||H|||_1=1, |||M|||_1^* \le 1\} < 1$ (4.2)

Let $X^*_S = \mathit{Arginf}_{\,\mathrm{supp}(S)\text{ in }S}|||Z|||_1 + (1/2)\gamma|\mathbf{y} - \Phi(Z)|_2^2$ be the minimizer of $MLP_{y,\Phi}(\gamma)$ with support in S, i.e.,

$$X^*_S = \mathit{Arg\ inf}\,|||Z|||_1 + (1/2)\gamma\,|\mathbf{y} - \Phi_S(Z)|_2^2 \tag{4.3}$$

Then there are the following conclusions:

(1) $X^*_S$ is the unique minimizer of problem (4.3) and is $l_1$-column-flat;

(2) $X^*_S$ is also the unique minimizer of problem (4.1), i.e., the (global) minimizer of (4.1) is unique and is $X^*_S$.

(3) Let $Y^* = \Phi_S^{*-1}(\mathbf{y}) \in \sum_s^{n\times n}(S)$, then $X^*_{S,ij} \ne 0$ for all $(i,j)$: $|Y^*_{ij}| > \gamma^{-1}N((\Phi_S^T\Phi_S)^{-1}: |||.|||_1^* \to |||.|||_{max})$ where $N((\Phi_S^T\Phi_S)^{-1}: |||.|||_1^* \to |||.|||_{max})$ denotes $(\Phi_S^T\Phi_S)^{-1}$'s operator norm and matrix norm $|||M|||_{max} := max_{ij}|M_{ij}|$.

(4) With the same notations as the above, if

$$min_{(i,j)\text{ in }S}|Y^*_{ij}| > \gamma^{-1}N((\Phi_S^T\Phi_S)^{-1}: |||.|||_1^* \to |||.|||_{max}) \tag{4.4}$$

then $sgn(X^*_{S,ij}) = sgn(Y^*_{ij})$ for all $(i,j)$ in S.

*Proof* : (1) Observe that when $\Phi_S^T\Phi_S$ is a bijection, (4.3)'s objective function $L_S(Z) = |||Z|||_1 + (1/2)\gamma\,|\mathbf{y} - \Phi_S(Z)|_2^2$ is strictly convex for variable $Z \in \sum_s^{n\times n}(S)$. According to general convex programming theory, its minimizer $X^*_S$ is unique.

(2) Let $L(Z) := |||Z|||_1 + (1/2)\gamma\,|\mathbf{y} - \Phi(Z)|_2^2$.To prove $X^*_S$ is also the global minimizer of (4.1), we prove its perturbation by H will always increase the objective's value, i.e., $L(X^*_S+H) > L(X^*_S)$ under the conditions specified by (1)(2)(3). Since conclusion (1) implies $L(X^*_S+H) > L(X^*_S)$ for any $H \ne O$ with support in S and $L(Z)$ is convex, we only need to consider the perturbation $X^*_S+H$ with $H_S = O$.

Since $X^*_S$ is the minimizer of (4.3), by first-order optimization condition there exists $M^*$ in $\partial|||X^*_S|||_1$ such that

$$M^* + \gamma\,\Phi_S^T(\Phi_S(X^*_S) - \mathbf{y}) = O \tag{4.5}$$

then $M^* = \gamma^*\Phi_S^T(\mathbf{y} - \Phi_S(X^*_S))$ and in particular $M^*_{\sim S} = O$. Equivalently:

$$X^*_S = \Phi_S^{*-1}(\mathbf{y}) - \gamma^{-1}(\Phi_S^T\Phi_S)^{-1}(M^*) \tag{4.6}$$

Now we compute $L(X^*_S+H) - L(X^*_S)$

$= |||X^*_S+H|||_1 - |||X^*_S|||_1 + (1/2)\gamma\,(|\Phi(X^*_S) - \mathbf{y}|_2^2 + 2<\Phi(X^*_S) - \mathbf{y}, \Phi(H)> + |\Phi(H)|_2^2 - |\Phi(X^*_S) - \mathbf{y}|_2^2)$

$= |||X^*_S+H|||_1 - |||X^*_S|||_1 + \gamma < \Phi(X^*_S) - \mathbf{y}, \Phi(H) > + (1/2)\gamma\,|\Phi(H)|_2^2$

$= |||X^*_S+H|||_1 - |||X^*_S|||_1 + \gamma < \Phi(X^*_S) - \mathbf{y}, \Phi_{\sim S}(H) > + (1/2)\gamma\,|\Phi_{\sim S}(H)|_2^2$

$\ge |||X^*_S+H|||_1 - |||X^*_S|||_1 + \gamma < \Phi(X^*_S) - \mathbf{y}, \Phi_{\sim S}(H) >$

The first term $|||X^*_S+H|||_1 - |||X^*_S|||_1$

$= max_j(|\mathbf{x}^*_j|_1 + |\mathbf{h}_j|_1) - max_j|\mathbf{x}^*_j|_1$ $\quad(\mathrm{supp}(X^*_S) \cap \mathrm{supp}(H) = \emptyset)$

$= |||X^*_S|||_1 + |||H|||_1 - |||X^*_S|||_1$ $\quad$ ( condition (1) implies $X^*_S$'s $l_1$-column-flatness: remark 3.5 )

$= |||\,H\,|||_1$

By replacing $X^*_S$ with (4.6), note $\mathrm{supp}(\Phi_{\sim S}{}^T) = \sim S$ and $|||M|||_1{}^* \le 1$, the second term

$\gamma <\Phi_S(X^*) - \mathbf{y}, \Phi_{\sim S}(H)>$

$= \gamma <\Phi_{\sim S}{}^T(\Phi_S\Phi_S{}^{*-1}(\mathbf{y}) - \mathbf{y}), H> - < M^*, \Phi_S{}^{*-1}\Phi_{\sim S}(H) >$ (4.7)

$\ge (-\gamma \sup\{<\Phi_{\sim S}{}^T(\Phi_S\Phi_S{}^{*-1}(\mathbf{y}) - \mathbf{y}), H>: |||H|||_1=1\} - \sup\{<\Phi_{\sim S}{}^T(\Phi_S{}^{*-1})^T M, H>: |||H|||_1=1, |||M|||_1{}^* \le 1\})|||H|||_1$

Therefore $L(X^*_S+H) - L(X^*_S) \ge |||\,H\,|||_1(1-\gamma \sup\{<\Phi_{\sim S}{}^T(\Phi_S\Phi_S{}^{*-1}(\mathbf{y}) - \mathbf{y}), H>: |||H|||_1 = 1\}$

$- \sup\{<\Phi_{\sim S}{}^T(\Phi_S{}^{*-1})^T M, H>: |||H|||_1 = 1 \text{ and } |||M|||_1{}^* \le 1\})$ (4.8)

and condition (3) implies the right hand side > 0. This proves $X^*_S$ is the minimizer of (4.1) and the minimizer is unique.

(3) For $Y^* = \Phi_S{}^{*-1}(\mathbf{y}) \in \Sigma_s^{n\times n}(S)$ (then $\mathrm{supp}(Y^*)$ in S) and by (4.6) we have

$|\,X^*_{S,ij}\,| = |\,Y^*_{ij} - \gamma^{-1}(\Phi_S{}^T\Phi_S)^{-1}(M^*)_{ij}\,| \ge |\,Y^*_{ij}\,| - \gamma^{-1}\,|\,(\Phi_S{}^T\Phi_S)^{-1}(M^*)_{ij}\,|$

$\ge |\,Y^*_{ij}\,| - \gamma^{-1}\, max_{ij}\,|\,(\Phi_S{}^T\Phi_S)^{-1}(M^*)_{ij}\,| = |\,Y^*_{ij}\,| - \gamma^{-1}\,|\,(\Phi_S{}^T\Phi_S)^{-1}(M^*)\,|_{max}$

$\ge |\,Y^*_{ij}\,| - N((\Phi_S{}^T\Phi_S)^{-1}: |||.|||_1{}^* \to |||.|||_{max})|||M|||_1{}^*$

$\ge |\,Y^*_{ij}\,| - N((\Phi_S{}^T\Phi_S)^{-1}: |||.|||_1{}^* \to |||.|||_{max})$ ( $|||M|||_1{}^* \le 1$ )

$> 0$ for those $(i,j)$: $|\,Y^*_{ij}\,| > N((\Phi_S{}^T\Phi_S)^{-1}: |||.|||_1{}^* \to |||.|||_{max})$

(4) Note that for any non-zero scalars $u$ and $v$, $sgn(u) = sgn(v)$ *iff* $|u| > |u - v|$. Therefore

$sgn(X^*_{S,ij}) = sgn(Y^*_{ij})$ *iff* $|\,Y^*_{ij}\,| > |\,Y^*_{ij} - X^*_{ij}\,| = \gamma^{-1}\,|\,(\Phi_S{}^T\Phi_S)^{-1}(M^*)_{ij}\,|$ (4.9)

In particular, if $min_{(i,j)\text{ in }S}|Y^*_{ij}| > \gamma^{-1}N((\Phi_S{}^T\Phi_S)^{-1}: |||.|||_1{}^* \to |||.|||_{max})$ then $|\,Y^*_{ij}\,| > \gamma^{-1}max_{(i,j)}\,|\,(\Phi_S{}^T\Phi_S)^{-1}(M^*)_{ij}\,|$ so $sgn(X^*_{S,ij}) = sgn(Y^*_{ij})$ for all $(i,j)$ in S. □

### 4.2 Conditions on Minimizer Uniqueness and Robustness for $MP_{y,\Phi,\eta}$

Now consider matrix signal reconstruction via solving the constrained convex programming:

$MP_{y,\Phi,\eta}$ : $X^* = Arg\ inf\,|||Z|||_1$ s.t. $Z \in R^{n\times n}$, $|\mathbf{y} - \Phi(Z)|_2 \le \eta$ (4.10)

**Lemma 4.2** Given **y**, positive integer $s$ and sparsity pattern $S = S_1\cup \dots \cup S_n$ where $|S_j| \le s$ for all $j$, suppose the linear measurement operator $\Phi$: $R^{n\times n} \to R^m$ satisfies:

(1) $\Phi_S{}^T(\mathbf{z})$ does not have any **0**-column for $\mathbf{z} \ne \mathbf{0}$;

(2) $\Phi_S{}^T\Phi_S$ is a bijection;

(3) $\sup\{<\Phi_{\sim S}{}^T(\Phi_S\Phi_S{}^{*-1}(\mathbf{y}) - \mathbf{y}), H>: |||H|||_1=1\} < \eta\, \Lambda_{min}(\Phi_S{}^T)\,(1 - N(\Phi_S{}^{*-1}\Phi_{\sim S}: |||.|||_1{}^* \to |||.|||_1{}^*))$

where $\Lambda_{min}(\Phi_S{}^T) := inf\,\{\,|||\Phi_S{}^T(z)|||_1{}^*: |z|_2 = 1\,\}$.

Then there are the following conclusions:

(1) As the minimizer of problem $MP_{y,\Phi,\eta}$, $X^*$ is unique, S-sparse and $l_1$-column-flat;

(2) Let $Y^* = \Phi_S{}^{*-1}(\mathbf{y}) \in \Sigma_s^{n\times n}(S)$, then for all $(i,j)$ in S:

$X^*_{ij} \ne 0$ and $sgn(X^*_{ij}) = sgn(Y^*_{ij})$ for all $(i,j)$: $|Y^*_{ij}| > \eta\, N(\Phi_S{}^T: l_2 \to |||.|||_1{}^*)N((\Phi_S{}^T\Phi_S)^{-1}: |||.|||_1{}^* \to |||.|||_{max})$

(3) For any given matrix norm $|.|_\alpha$ there holds:

$$|X^* - Y^*|_\alpha \le \eta\, N(\Phi_S{}^{*-1}: l_2 \to |.|_\alpha) \quad (4.11)$$

*Proof* (1) Let $X^*_S \in Arg\ inf\,|||Z|||_1$ s.t. $Z \in R^{n\times n}$, $|\mathbf{y} - \Phi_S(Z)|_2 \le \eta$. i.e., a minimizer with its support restricted on S. We first prove $X^*_S$ is the only minimizer of this support-restricted problem, then we prove $X^*$ is also the minimizer of problem $MP_{y,\Phi,\eta}$ (4.10),i.e., $X^*_S$ is the global minimizer and (4.10)'s minimizer is unique.

According to general convex optimization theory, there exist a positive multiplier $\gamma^* > 0$ and $M^*$ in $\partial|||X^*_S|||_1$ such that $M^* + \gamma^*\Phi_S{}^T(\Phi_S(X^*_S) - \mathbf{y}) = O$ and $|\mathbf{y} - \Phi_S(X^*_S)|_2 = \eta$ (4.12)

then equivalently $X^*_S = \Phi_S{}^{*-1}(\mathbf{y}) - \gamma^{*-1}(\Phi_S{}^T\Phi_S)^{-1}(M^*)$ (4.13)

Suppose $X^0$ is another minimizer of $inf\,|||Z|||_1$ s.t. $Z \in R^{n\times n}$, $|\mathbf{y} - \Phi_S(Z)|_2 \le \eta$, then there exist a positive multiplier $\gamma^0 > 0$ and $M^0$ in $\partial|||X^0|||_1$ such that

$$M^0 + \gamma^0\Phi_S^T(\Phi_S(X^0) - \mathbf{y}) = O \text{ and } |\mathbf{y} - \Phi_S(X^0)|_2 = \eta \tag{4.14}$$

Equivalently, (4.12) shows that $X^*_S$ is also a minimizer of $L_S(Z) = |||Z|||_1 + (1/2)\gamma^*|\mathbf{y} - \Phi_S(Z)|_2^2$ which is a strictly convex function on $\sum_S^{n\times n}(S)$ since $\Phi_S^T\Phi_S$ is a bijection (condition(2)), as a result $L_S(Z)$'s minimizer is unique. However, since $|||X^*_S|||_1 = |||X^0|||_1$ we have $L_S(X^*_S) = |||X^*_S|||_1 + (1/2)\gamma^*|\mathbf{y} - \Phi_S(X^*_S)|_2^2 = |||X^*_S|||_1 + \gamma^*\eta^2/2 = |||X^0|||_1 + (\gamma^*/2)|\mathbf{y} - \Phi_S(X^0)|_2^2 = L_S(X^0)$, which implies $X^*_S = X^0$, i.e., $X^*_S$ is the unique minimizer of the support-restricted problem $inf\,|||Z|||_1$ s.t. $Z \in R^{n\times n}$, $|\mathbf{y}-\Phi_S(Z)|_2 \le \eta$.

$X_S^*$'s $l_1$-column-flatness is implied by condition (1) and theorem 3.3.

Now prove $X^*_S$ (which is S-sparse and $l_1$-column-flat ) is also a minimizer of problem $MP_{y,\Phi,\eta}$ (4.10). Again we start with the fact that $X^*_S = Arginf\ L_S(Z) = Arginf\ |||Z|||_1 + (1/2)\gamma^*\ |\mathbf{y} - \Phi_S(Z)|_2^2$ with some multiplier $\gamma^* > 0$ (which value depends on $X^*_S$) and by lemma 4.1, $X^*_S$ is the unique minimizer of the convex problem (without any restriction on solution's support )

$$inf\,|||Z|||_1 + (1/2)\gamma^*\ |\mathbf{y} - \Phi(Z)|_2^2 \tag{4.15}$$

under the condition

$$\gamma^* \sup\{<\Phi_{\sim S}^T(\Phi_S\Phi_S^{*-1}(\mathbf{y}) - \mathbf{y}), H>: |||H|||_1=1\}$$
$$+ \sup\{<\Phi_{\sim S}^T(\Phi_S^{*-1})^T(M),H>: |||H|||_1=1 \text{ and } |||M|||_1^* \le 1\} < 1 \tag{4.16}$$

According to convex optimization theory, $X^*_S$ (under condition (4.16)) being the unique minimizer of problem (4.15) means $X^*_S$ is also a minimizer of $MP_{y,\Phi,\eta}$ (4.10), which furthermore implies that $MP_{y,\Phi,\eta}$'s minimizer is unique, S-sparse and $l_1$-column-flat.

In order to make condition (4.16) more meaningful, we need to replace the minimizer-dependent parameter $\gamma^*$ with explicit information. From (4.15)'s first-order optimization condition (4.12) we obtain

$$1 \ge |||M^*|||^*_1 = \gamma^*\ |||\Phi_S^T(\Phi_S(X^*_S) - \mathbf{y})\ |||^*_1 \ge \gamma^*\ min\{\ |||\Phi_S^T(\mathbf{z})|||^*_1: |\mathbf{z}|_2 = 1\}|\Phi_S(X^*_S) - \mathbf{y}|_2 = \gamma^*\eta\ \Lambda_{min}(\Phi_S^T)$$

i.e.,
$$\gamma^* \le (\eta\ \Lambda_{min}(\Phi_S^T))^{-1} \tag{4.17}$$

with this upper-bound of $\gamma^*$, (4.16) can be derived from a uniform condition

$$(\eta\ \Lambda_{min}(\Phi_S^T))^{-1} \sup\{<\Phi_{\sim S}^T(\Phi_S\Phi_S^{*-1}(\mathbf{y}) - \mathbf{y}), H>: |||H|||_1=1\}$$
$$+ \sup\{<\Phi_{\sim S}^T(\Phi_S^{*-1})^T(M), H>: |||H|||_1=1 \text{and } |||M|||_1^* \le 1\} < 1 \tag{4.18}$$

which is equivalent to condition (3).

From now on we denote $X^*_S$ as $X^*$.

(2) For $Y^* = \Phi_S^{*-1}(\mathbf{y}) \in \sum_S^{n\times n}(S)$ and by lemma 4.1's conclusion (4), if $min_{(i,j)\text{ in }S}|Y^*_{ij}| > \gamma^{*-1}N((\Phi_S^T\Phi_S)^{-1}: |||.|||_1^* \to |||.|||_{max})$ then $sgn(X^*_{S,ij}) = sgn(Y^*_{ij})$ for all $(i,j)$ in S. To replace multiplier $\gamma^*$ with more explicit information in this condition, we need some lower bound of $\gamma^*$ which can be derived from the first-order optimization condition $M^* = \gamma^*(\mathbf{y} - \Phi_S^T(\Phi_S(X^*))$ again. Note that $X^*$ is $l_1$-column-flat implies every column of $X^*$ is not **0**, further more $M^*$ has no **0**-column so $M^* = (\lambda_1\mathbf{u}_1,\dots, \lambda_n\mathbf{u}_n)$ with $\lambda_j > 0$ for all $j$, $\lambda_1+\dots+\lambda_n = 1$ and $|\mathbf{u}_j|_\infty = 1$, as a result $|||M^*|||_1^* = \sum_j\lambda_j|\mathbf{u}_j|_\infty = 1$. Hence

$$1 = |||M^*|||_1^* \le \gamma^*\ |||\Phi_S^T(\Phi_S(X^*) - \mathbf{y})\ |||_1^* \le \gamma^*\ N(\Phi_S^T: l_2 \to |||.|||_1^*)|\Phi_S(X^*) - \mathbf{y}|_2 = \gamma^*\eta\ N(\Phi_S^T: l_2 \to |||.|||_1^*)$$

i.e.,
$$\gamma^{*-1} \le \eta\ N(\Phi_S^T: l_2 \to |||.|||_1^*) \tag{4.19}$$

Replace $\gamma^{*-1}$ with its upper-bound in (4.19), we obtain if $min_{(i,j)\text{ in }S}\ |Y^*_{ij}| > \eta\ N(\Phi_S^T: l_2 \to |||.|||^*_1)N((\Phi_S^T\Phi_S)^{-1}: |||.|||_1^* \to |||.|||_{max})$ then $sgn(X^*_{S,ij}) = sgn(Y^*_{ij})$ for all $(i,j)$ in S.

(3) $Y^* = \Phi_S^{*-1}(\mathbf{y}) \in \sum_S^{n\times n}(S)$ implies $\Phi_S^T(\Phi_S(Y^*) - \mathbf{y}) = O$ and then condition (1) leads to $\Phi_S(Y^*) = \mathbf{y}$. Furthermore, $\Phi_S^T\Phi_S$ is a bijection for $\sum_S^{n\times n}(S) \to \sum_S^{n\times n}(S)$ and notice $X^* - Y^* \in \sum_S^{n\times n}(S)$, so for any matrix norm $|.|_\alpha$:

$$|X^* - Y^*|_\alpha = |(\Phi_S^T\Phi_S)^{-1}(\Phi_S^T\Phi_S)(X^* - Y^*)|_\alpha = |\Phi_S^{*-1}\Phi_S(X^* - Y^*)|_\alpha = |(\Phi_S^{*-1}(\Phi_S(X^*) - \mathbf{y}))|_\alpha$$
$$\le N(\Phi_S^{*-1}: l_2 \to |.|_\alpha)|\Phi_S(X^*) - \mathbf{y}|_2 = \eta\ N(\Phi_S^{*-1}: l_2 \to |.|_\alpha)$$ □

*Remark* 4.1: Under the conditions specified in this lemma, the minimizer $X^*$ can have support stability, sign stability and component value-error robustness relative to any metric $|.|_\alpha$, e.g., we can take $|.|_\alpha$ as $|.|_F$, $|||.|||_1$, $|||.|||_{max}$, etc. The error is linearly upper-bounded by the measurement error $\eta$. When $\eta$ decreases as small as possible, the error $|X^* - Y^*|_\alpha$ seems to be reduced also as small as possible. However, this is not true in general because condition (3) may require $\eta$ not be two small.

In real world applications, the support S of the signal is of course unknown so lemma 4.2 cannot be applied directly. However, on basis of lemma 4.2 a stronger and uniform sufficient condition can be established to guarantee the uniqueness and robustness of the reconstructed signal from solving $MP_{y,\Phi,\eta}$.

**Theorem 4.1** Given positive integer $s$ and the linear measurement operator $\Phi$: $R^{n\times n} \rightarrow R^m$, suppose $\Phi$ satisfies the following conditions for any $s$-sparsity pattern $S = S_1 \cup \ldots \cup S_n$ where $|S_j| \le s$ for all $j$:

(1) $\Phi_S^T(\mathbf{z})$ does not have any **0**-column for $\mathbf{z} \neq \mathbf{0}$;

(2) $\Phi_S^T\Phi_S$ is a bijection;

(3) $N(\Phi_{\sim S}^T(\Phi_S\Phi_S^{*-1} - I_S): l_2 \rightarrow |||.|||_1^*) < \Lambda_{min}(\Phi_S^T)\,(1 - N(\Phi_S^{*-1}\Phi_{\sim S}: |||.|||_1^* \rightarrow |||.|||_1^*))$ or equivalently

$N(((\Phi_S^{*-1})^T\Phi_S^T - I_S^T)\Phi_{\sim S} : |||.|||_1 \rightarrow l_2) < \Lambda_{min}(\Phi_S^T)\,(1 - N(\Phi_{\sim S}^T(\Phi_S^{*-1})^T: |||.|||_1 \rightarrow |||.|||_1))$

where $\Lambda_{min}(\Phi_S^T) := inf\{\, |||\Phi_S^T(z)|||_1^* : |z|_2 = 1 \,\}$ and $I_S :=$ the identical mapping on $\sum_s^{n\times n}(S)$.

Then for the minimizer $X^*$ of problem $MP_{y,\Phi,\eta}$ (4.10) where $\mathbf{y} = \Phi(X) + \boldsymbol{e}$ with noise $|\boldsymbol{e}|_2 \le \eta$ and a real flat signal $X \in \sum_s^{n\times n}(R)$ of some $s$-sparsity pattern R, there are the following conclusions:

(1) *Sparsity, flatness and support stability*:

$X^* \in \sum_s^{n\times n}(R)$ and is $l_1$-column-flat and the unique minimizer of $MP_{y,\Phi,\eta}$;

(2) *Robustness*: For any given matrix norm $|.|_\alpha$ there holds:

$$|X^* - X|_\alpha \le 2\eta\, N(\Phi_R^{*-1}: l_2 \rightarrow |.|_\alpha) \qquad (4.20)$$

(3) *Sign Stability*: $sgn(X^*_{ij}) = sgn(X_{ij})$ for $(i,j)$ in R such that :

$$|X_{ij}| > \eta\,(\, N(\Phi_R^{*-1}: l_2 \rightarrow |||.|||_{max}) + N(\Phi_R^T: l_2 \rightarrow |||.|||_1^*) N((\Phi_R^T\Phi_R)^{-1}: |||.|||_1^* \rightarrow |||.|||_{max})\,) \qquad (4.21)$$

*Proof* (1) Note that in case of $X \in \sum_s^{n\times n}(R)$ and $\mathbf{y} = \Phi(X) + \boldsymbol{e} = \Phi_R(X) + \boldsymbol{e}$, $|\boldsymbol{e}|_2 \le \eta$, we have

$$\Phi_R\Phi_R^{*-1}(\mathbf{y}) - \mathbf{y} = (\Phi_R\Phi_R^{*-1} - I_R)\boldsymbol{e}$$

It's straightforward to verify that in this situation condition (3) in this theorem leads to condition (3) in lemma 4.2: $\sup\{\langle \Phi_{\sim R}^T(\Phi_R\Phi_R^{*-1}(\mathbf{y}) - \mathbf{y}), H\rangle : |||H|||_1 = 1\} < \eta\, \Lambda_{min}(\Phi_R^T)\,(1 - N(\Phi_R^{*-1}\Phi_{\sim R}: |||.|||_1^* \rightarrow |||.|||_1^*))$ for any $\eta$. As a result, $X^* \in \sum_s^{n\times n}(R)$ and is $l_1$-column-flat and the unique minimizer of $MP_{y,\Phi,\eta}$.

(2) For $Y^* = \Phi_R^{*-1}(\mathbf{y}) \in \sum_s^{n\times n}(R)$ and by lemma 4.2(4), we obtain $|X^* - Y^*|_\alpha \le \eta\, N(\Phi_R^{*-1}: l_2 \rightarrow |.|_\alpha)$ for any given matrix norm $|.|_\alpha$. On the other hand, $Y^* = \Phi_R^{*-1}(\mathbf{y})$ implies $\Phi_R^T(\Phi_R(Y^*) - \mathbf{y}) = O$ then condition (1) leads to $\Phi_R(Y^*) = \mathbf{y}$, hence $\Phi_R(Y^*) = \mathbf{y} = \Phi(X) + \boldsymbol{e} = \Phi_R(X) + \boldsymbol{e}$, namely $\Phi_R^T\Phi_R(Y^*) = \Phi_R^T\Phi_R(X) + \Phi_R^T(\boldsymbol{e})$, as a result:

$$Y^* - X = (\Phi_R^T\Phi_R)^{-1}\Phi_R^T(\boldsymbol{e}) \equiv \Phi_R^{*-1}(\boldsymbol{e}) \qquad (4.22)$$

Since $|\boldsymbol{e}|_2 \le \eta$, we get $|Y^* - X|_\alpha \le \eta\, N(\Phi_R^{*-1}: l_2 \rightarrow |.|_\alpha)$ for any given matrix norm $|.|_\alpha$. Combining with $|X^* - Y^*|_\alpha \le \eta\, N(\Phi_R^{*-1}: l_2 \rightarrow |.|_\alpha)$ we get the reconstruction error bound $|X^* - X|_\alpha \le 2\eta\, N(\Phi_R^{*-1}: l_2 \rightarrow |.|_\alpha)$

(3) By the first-order optimization condition on minimizer $X^*$ with the fact $supp(X^*) = R$, we have the equation $X^* = \Phi_R^{*-1}(\mathbf{y}) - \gamma^{*-1}(\Phi_R^T\Phi_R)^{-1}(M^*) = Y^* - \gamma^{*-1}(\Phi_R^T\Phi_R)^{-1}(M^*)$ where $M^*$ is in $\partial|||X^*|||_1$, namely:

$$X^* - Y^* = -\gamma^{*-1}(\Phi_R^T\Phi_R)^{-1}(M^*) \qquad (4.23)$$

Combining with (4.22), we get

$$X^* - X = \Phi_R^{*-1}(\boldsymbol{e}) - \gamma^{*-1}(\Phi_R^T\Phi_R)^{-1}(M^*) \qquad (4.24)$$

Since $sgn(X^*_{ij}) = sgn(X_{ij})$ *iff* $|X_{ij}| > |X_{ij} - X^*_{ij}| = |\Phi_R^{*-1}(\boldsymbol{e})_{ij} - \gamma^{*-1}(\Phi_R^T\Phi_R)^{-1}(M^*)_{ij}|$, in particular, if $X_{ij}$ can satisfy $|X_{ij}| > max_{ij}|\Phi_R^{*-1}(\boldsymbol{e})_{ij}| + \gamma^{*-1} max_{ij}|(\Phi_R^T\Phi_R)^{-1}(M^*)_{ij}|$ then the former inequality is true and as a result $sgn(X^*_{ij}) = sgn(X_{ij})$. It's straightforward to verify (by using (4.19)) that the condition (4.21) just provides a guarantee for this. □

*Remark* 4.2: In this theorem condition (3) is independent with measurement error bound $\eta$, which implies $|X^* - X|_\alpha = O(\eta)$ can hold for any small value of $\eta$, and (4.21) indicates that with $\eta$ small enough, all signal's

non-zero components' signs can be correctly recovered. For given Φ, the associated value

$$Max_{s\text{-sparse pattern S}}\, N(\Phi_S^{*-1}: l_2 \to |||.|||_{max}) + N(\Phi_S^T: l_2 \to |||.|||_1^*) N((\Phi_S^T\Phi_S)^{-1}: |||.|||_1^* \to |||.|||_{max})$$

or an enhancement $2Max_{s\text{-sparse pattern S}}\, N(\Phi_S^T: l_2 \to |||.|||_1^*) N((\Phi_S^T\Phi_S)^{-1}: |||.|||_1^* \to |||.|||_{max})$ can be regarded as a signal-to-noise ratio threshold for correct sign recovery.

Finally we mention that in condition (3) $\Lambda_{min}(\Phi_S^T)=inf\{|||\Phi_S^T(z)|||_1^*: |z|_2=1\}$ for $MP^{(2)}_{y,\Phi,\eta}$ or $\Lambda_{min}(\Phi_S^T)=inf\{|||\Phi_S^T(Z)|||_1^*: |Z|_F=1\}$ for $MP^{(F)}_{Y,A,B,\eta}$ is one of critical quantities which value should be as large as possible to satisfy this condition. Observe that there is a counterpart $\lambda_{min,\alpha\beta}(\Phi;K)=inf\{|\Phi(Z)|_\alpha: Z \text{ in } K \text{ and } |Z|_\beta=1\}$ (e.g., see FACT2.6) which is also critical for analysis in next sections in random measurement setting.

# 5 Number of Measurements for Robust Reconstruction via Solving $MP^{(2)}_{y,\Phi,\eta}$

After established the conditions for the measurement operator to guarantee desired properties (e.g., uniqueness, robustness, etc.) of the matrix signal reconstruction, next question should be about how to construct the measurement operator to satisfy such conditions with required number of measurements as few as possible. This and next sections deal with this problem in random approach. In this section we establish conditions on number of measurements $m$ for robustly reconstructing the matrix signal X by solving the convex programming problem $MP^{(2)}_{y,\Phi,\eta}$:

$$inf\,|||Z|||_1 \quad \text{s.t } Z \in R^{n\times n},\ |\mathbf{y} - \Phi(Z)|_2 \le \eta$$

where $\mathbf{y} = \Phi(X) + \boldsymbol{e}$ and $|\boldsymbol{e}|_2 \le \eta$, $\Phi$: $R^{n\times n} \to R^m$ is a linear operator. In the equivalent component-wise formulation, $y_i = \langle\Phi_i, X\rangle + e_i$ where for $i=1,\ldots,m$, $\Phi_i \in R^{n\times n}$ are random matrices independent each other and $\Phi_i$'s entries are independently sampled under standard Gaussian $N(0,1)$ or sub-Gaussian distribution. In section 5.3, a simple necessary condition on $m$ is established.

Instead of proving that M-RIP property (defined in sec. 3.5) can be satisfied with high probability in context of the distributions in consideration, which are quite involved particularly for problem $\mathbf{MP^{(F)}_{Y,A,B,\eta}}$, here we straightforwardly prove that the robustness in terms of *Frobenius*-norm error metric can be reached with high probability when number of measurements exceeds some specific bound. The investigation on how those conditions established in last sections can be satisfied by the random measurement operator will be subjects in subsequent papers.

## 5.1 Case 1: Gaussian Measurement Operator Φ

Based upon the fundamental facts presented in section 2.2, one of the critical steps in this approach is to estimate the width $w(D(|||X|||_1, X))$'s upper bound for matrix signal $X=(\mathbf{x}_1,\ldots,\mathbf{x}_n)$ with $s$-sparse column vectors $\mathbf{x}_1,\ldots,\mathbf{x}_n$. This is done in lemma 5.1.

Based on lemma 3.1 and FACT 2.4, the upper bound of Gaussian width $D(|||.|||_1, X)$ with respect to *Frobenius* norm is estimated in the following lemma.

**Lemma 5.1** Given $n$-by-$n$ matrix $X=(\mathbf{x}_1,\ldots,\mathbf{x}_n)$ with $s$-sparse column vectors $\mathbf{x}_1,\ldots,\mathbf{x}_n$. Let $r$ (called *$l_1$-column-flatness parameter* hereafter) be cardinality of the set $\{j: |\mathbf{x}_j|_1=max_k|\mathbf{x}_k|_1\}$, i.e., the number of column vectors which have the maximum $l_1$-norm. Then

$$w^2(D(|||.|||_1, X)) \le 1+n^2-r(n-slog(Cn^4r^2)) \qquad (5.1)$$

In particular, when $r = n$ then

$$w^2(D(|||.|||_1, X)) \le 1 + nslog(Cn^6) \qquad (5.2)$$

where $C$ is an absolute constant.

*Remark* 5.1: $ns$ is the total sparsity of the matrix signal X. This estimate shows that the signal complexity

(width) encoded by regularizer $|||X|||_1$ is controlled by two structural parameters, the column-sparsity $s$ and $l_1$-column-flatness $r$. Complexity gets lower with smaller $s$ and larger $r$.

*Proof* We start with (FACT 2.4) $w^2(D(|||.|||_1, X)) \le E_G[inf\{|G-tV|_F^2: t>0, V \text{ in } \partial|||X|||_1\}]$ where G is a random matrix with entries $G_{ij} \sim^{iid} N(0,1)$.

Set $G=(\mathbf{g}_1,\ldots,\mathbf{g}_n)$ where $\mathbf{g}_j \sim^{iid} N(0,I_n)$. By lemma 3.1, $V=(\lambda_1\xi_1,\ldots, \lambda_n\xi_n)$ where w.l.o.g. $\lambda_j \ge 0$ for $j=1,\ldots,r$, $\lambda_1+\ldots+\lambda_r=1$, $\lambda_j=0$ for $j\ge r+1$; $|\mathbf{x}_j|_1=max_k|\mathbf{x}_k|_1$ for $j=1,\ldots,r$ and $|\mathbf{x}_j|_1<max_k|\mathbf{x}_k|_1$ for $j\ge 1+r$; $\xi_j(i)=sgn(X_{ij})$ for $X_{ij}\ne 0$ and $|\xi_j(i)| \le 1$ for all $i$ and $j$. Then

$$w^2(D(|||.|||_1, X)) \le E_G[inf_{t>0,\ \lambda j,\ \xi j \text{ specified as the above}} \textstyle\sum_{j=1}^{r} |\mathbf{g}_j - t\lambda_j\xi_j|_2^2 + \sum_{j=r+1}^{n} |\mathbf{g}_j|_2^2]$$

$$\le inf_{t>0, \text{ all } \lambda j \text{ specified as the above}} E_G[inf_{\text{all } \xi j \text{ specified as the above}} \textstyle\sum_{j=1}^{r} |\mathbf{g}_j - t\lambda_j\xi_j|_2^2 + \sum_{j=r+1}^{n} |\mathbf{g}_j|_2^2]$$

$$= inf_{t>0, \text{ all } \lambda j \text{ specified as the above}} E_G[inf_{\text{all } \xi j \text{ specified as the above}} \textstyle\sum_{j=1}^{r} |\mathbf{g}_j - t\lambda_j\xi_j|_2^2]+\sum_{j=r+1}^{n} E_G[|\mathbf{g}_j|_2^2]$$

$$= inf_{t>0, \text{ all } \lambda j \text{ specified as the above}} E_G[\textstyle\sum_{j=1}^{r} inf_{\xi j \text{ specified as the above}} |\mathbf{g}_j - t\lambda_j\xi_j|_2^2] + (n-r)n$$

( since $\xi_j$ is unrelated each other and $E_G[|\mathbf{g}_j|_2^2]=n$ )

$$= inf_{t>0, \text{ all } \lambda j \text{ specified as the above}} \textstyle\sum_{j=1}^{r} E_{gj}[inf_{\xi j \text{ specified as the above}} |\mathbf{g}_j - t\lambda_j\xi_j|_2^2] + (n-r)n$$

For each $j=1,\ldots,r$ let $S(j)$ be the support of $\mathbf{x}_j$ (so $|S(j)| \le s$) and $\sim S(j)$ be its complimentary set, then $|\mathbf{g}_j - t\lambda_j\xi_j|_2^2=|\mathbf{g}_{j|S(j)} - t\lambda_j\xi_{j|S(j)}|_2^2 + |\mathbf{g}_{j|\sim S(j)} - t\lambda_j\xi_{j|\sim S(j)}|_2^2$. Notice that all components of $\xi_{j|S(j)}$ are $\pm 1$ and all components of $\xi_{j|\sim S(j)}$ can be any value in the interval $[-1,+1]$. Select $\lambda_1=\ldots=\lambda_r=1/r$, let $\varepsilon>0$ be arbitrarily small positive number and select $t=t(\varepsilon)$ such that $P[|g|>t(\varepsilon)/r]\le\varepsilon$ where $g$ is a standard scalar Gaussian random variable (i.e., $g\sim N(0,1)$ and $\varepsilon$ can be $exp(-t(\varepsilon)^2/2r^2)$). For each $j$ and each $i$ outside $S(j)$, set $\xi_j$'s component $\xi_j(i) = rg_j(i)/t(\varepsilon)$ if $|g_j(i)| \le t(\varepsilon)/r$ (in this case $|g_j(i) - t\lambda_j\xi_{j|}(i)| = 0$) and otherwise $\xi_j(i)=sgn(g_j(i))$ (in this case $|g_j(i) - t\lambda_j\xi_{j|}(i)| = |g_j(i)| - t(\varepsilon)/r$), then $|\mathbf{g}_{j|\sim S(j)} - t\lambda_j\xi_{j|\sim S(j)}|_2^2=0$ when $|\mathbf{g}_{j|\sim S(j)}|_\infty < t(\varepsilon)/r$, hence:

$$E[ |\mathbf{g}_{j|\sim S(j)} - t\lambda_j\xi_{j|\sim S(j)}|_2^2] = \int_0^\infty du P[ |\mathbf{g}_{j|\sim S(j)} - t\lambda_j\xi_{j|\sim S(j)}|_2^2>u] = 2\int_0^\infty duuP[ |\mathbf{g}_{j|\sim S(j)} - t\lambda_j\xi_{j|\sim S(j)}|_2>u]$$

$$\le 2 \int_0^\infty duuP[\text{There exists } (\mathbf{g}_{j|\sim S(j)} - t\lambda_j\xi_{j|\sim S(j)})\text{'s component with magnitude} > (n-s)^{-1/2}u ]$$

$$\le 2(n-s) \int_0^\infty duuP[ |g| - t(\varepsilon)/r > (n-s)^{-1/2}u ]$$

$$\le 2(n-s) \int_0^\infty duuexp(-((t(\varepsilon)/r)+ (n-s)^{-1/2}u)^2/2)$$

$$\le C_0(n-s)^2exp(-t(\varepsilon)^2/2r^2 ) \le C_0(n-s)^2\varepsilon$$

where $C_0$ is an absolute constant. On the other hand:

$$E_{gj}[|\mathbf{g}_{j|S(j)} - t\lambda_j\xi_{j|S(j)}|_2^2] = E_{gj}[|\mathbf{g}_{j|S(j)}|^2] + (t(\varepsilon)^2/r^2)|\xi_{j|S(j)}|_2^2 = (1+ t(\varepsilon)^2/r^2)s=(1+2log\ (1/\varepsilon))s$$

Hence $w^2(D(|||.|||_1, X)) \le (1+2log\ (1/\varepsilon))rs + (n-r)n + r(n-s)^2\varepsilon \le n^2 - r(n-slog(e/\varepsilon^2)) + C_0n^2r\varepsilon$

In particular, let $\varepsilon=1/C_0n^2r$ then we get $w^2(D(|||.|||_1, X)) \le n^2 - r(n - slog(Cn^4r^2)) + 1$. □

Combing this lemma and FACT 2.3, we obtain the general result in the following:

**Theorem 5.1** Suppose $\Phi_{kij} \sim^{iid} N(0,1)$, let $X\in \sum_s^{n\times n}$ be a columnwise $s$-sparse and $l_1$-column-flat matrix signal, $R^m \ni \mathbf{y} = \Phi(X) + e$ where $|e|_2 \le \eta$, $X^*$ be the minimizer of the problem $MP^{(2)}_{y,\Phi,\eta}$. If the measurement vector $\mathbf{y}$'s dimension

$$m \ge (t + 2\eta/\delta + (nslog(Cn^6))^{1/2})^2$$

where $C$ is an absolute constant, then $P[|X^*-X|_F \le \delta ] \ge 1 - exp(-t^2/2)$, i.e., such X can be reconstructed robustly with respect to the error norm $|X^*-X|_F$ with high probability by solving $MP^{(2)}_{y,\Phi,\eta}$. □

## 5.2 Case 2: Sub-Gaussian Measurement Operator Φ

Combing lemma 5.1, FACT 2.1(2) and 2.5, the following result can be obtained straightforwardly:

**Theorem 5.2** Let X and $X^*$ be respectively the matrix signal and the minimizer of $MP^{(2)}_{y,\Phi,\eta}$ where the signal $X\in \sum_s^{n\times n}$ is columnwise $s$-sparse and $l_1$-column-flat, $y_k=\langle\Phi_k,X\rangle+e_k$, each $\Phi_k \sim^{iid} \Phi$ where $\Phi$ is a random matrix satisfying the conditions (1)(2)(3) in FACT 2.5 with parameters $\alpha$, $\rho$ and $\sigma$. If the measurement vector **y**'s dimension

$$m \geq (C_1\rho^4/\alpha)(\ \alpha t + 2\eta/\delta + \sigma C_2(\rho^6 nslog(C_3n^6))^{1/2})^2$$

where $C_i$'s are absolute constants, then $P[|X^*-X|_F \leq \delta] \geq 1 - exp(-C_4t^2)$, i.e., such X can be reconstructed robustly with respect to the error norm $|X^*-X|_F$ with high probability by solving $MP^{(2)}_{y,\Phi,\eta}$. □

### 5.3 Necessary Condition on Number of Measurements

**Theorem 5.3** Given measurement operator $\Phi$: $R^{n\times n} \to R^m$, if any $X=(\mathbf{x}_1,\ldots,\mathbf{x}_n)$ with sparse column vectors $\mathbf{x}_j$ in $\sum^{2S}$ for all $j$ is always the unique solution to problem $MP^{(\alpha)}_{y,\Phi,\eta}$ where $\eta=0$, then

$$m \geq C_1 nslog(C_2n/s)) \quad (5.3)$$

where $C_1$ and $C_2$ are absolute constants.

*Proof* For any $s<n$, there exist $k \geq (n/4s)^{ns/2}$ subsets $S^{(\alpha\beta\ldots\omega)}=S_1^{(\alpha)}\cup S_2^{(\beta)}\cup\ldots\cup S_n^{(\omega)}$ in $\{(i,j): 1\leq i, j\leq n\}$ where each $S_j^{(\mu)}=\{(i_1,j),\ldots,(i_s,j): 1\leq i_1<i_2<\ldots<i_s\leq n\}$ and $|S_j^{(\mu)}\cap S_j^{(\nu)}| < s/2$ for $\mu\neq\nu$. This fact is based on a combinatorial theorem[11] that for any $s<n$ there exist $l \geq (n/4s)^{s/2}$ subsets $R^{(\mu)}$ in $\{1,2,\ldots,n\}$ where $|R^{(\mu)}\cap R^{(\nu)}| < s/2$ for any $\mu\neq\nu$. For the $n$-by-$n$ square $\{(i,j): 1\leq i, j\leq n\}$, assign a $R^{(\mu)}$ to each column, i.e., set $S_j^{(\mu)}:=\{(i, j): i\in R^{(\mu)}\}$. As a result $|S_j^{(\mu)}\cap S_j^{(\nu)}| < s/2$ for $\mu\neq\nu$ since $|R^{(\mu)}\cap R^{(\nu)}| < s/2$ for $\mu\neq\nu$ and totally there can be $k=l^n$ such assignments $S^{(\alpha\beta\ldots\omega)}=S_1^{(\alpha)}\cup S_2^{(\beta)}\cup\ldots\cup S_n^{(\omega)}$ on the square.

Now we call the above $S_1^{(\alpha)}\cup S_2^{(\beta)}\cup\ldots\cup S_n^{(\omega)}$ a *configuration* on the $n$-by-$n$ square. Let $m$ be the rank of linear operator $\Phi$. Consider the quotient space $L:=R^{n\times n}/\ker\Phi$, then $\dim L=n^2-\dim\ker\Phi=m$. For any [X] in $L$ define the norm $|[X]|:=inf\{|||X-V|||_1: V \text{ in } \ker\Phi\}$. For any $X=(\mathbf{x}_1,\ldots,\mathbf{x}_n)$ with $\mathbf{x}_j$ in $\sum^{2S}$ for all $j$ , the assumption about $\Phi$ implies $|[X]|=|||X|||_1$. Now for any configuration $\Delta=S_1\cup S_2\cup\ldots\cup S_n$ on the $n$-by-$n$ square, define $X_{ij}(\Delta):=1/s$ if $(i,j)\in S_j$ and 0 otherwise, then $|||X(\Delta)|||_1=1$, each $X(\Delta)$'s column $x_j(\Delta)\in\sum^S$ and each column of $X(\Delta')-X(\Delta'')$ is in $\sum^{2S}$, furthermore $|[X(\Delta')]-[X(\Delta'')]|=|||X(\Delta')-X(\Delta'')|||_1>1$ because of the property $|S_j'\cap S_j''| < s/2$ for $S_j'\neq S_j''$. These facts imply that the set $\Theta:=\{[X(\Delta)]: \Delta \text{ runs over all configurations}\}$ is a subset on normed quotient space $L$'s unit sphere with distances between any pair of its members >1, i.e., a $d$-net on the sphere where $d>1$. The cardinality of $\Theta$ = number of configurations $k \geq (n/4s)^{ns/2}$ and an elementary estimate derives $k \leq 3^{dimL}=3^m$, hence $m \geq C_1nslog(C_2n/s))$ where $C_1=1/2log3$ and $C_2=1/4$. □

*Remark* 5.2: For the measurement operator $\Phi_{A,B}$: $R^{n\times n} \to R^{m\times m}$:$Y=AXB^T$, the same result is true with $m^2 \geq C_1nslog(C_2n/s))$.

## 6 Number of Measurements for Robust Reconstruction via Solving $MP^{(F)}_{Y, A, B, \eta}$

Now we investigate the problem $MP^{(F)}_{Y, A, B, \eta}$ in which A and B are both Gaussian or sub-Gaussian $m$-by-$n$ random matrices. Notice that in these cases the measurement operator $\Phi_{A,B}$ is neither Gaussian nor sub-Gaussian. The basis is FACT 2.6 and the critical step is also the width's upper bound estimation.

More explicitly, $\Phi_{A,B}$: $R^{n\times n} \to R^{m\times m}$: $Y=AZB^T$. In the equivalent component-wise formulation, $y_{kl}=\langle\Phi_{kl},X\rangle=\sum_{ij}A_{ki}X_{ij}B_{lj}$ for each $1\leq k,l\leq m$, $A_{ki} \sim^{iid} B_{lj} \sim^{iid} N(0,1)$ or sub-Gaussian distribution, and A, B are independent each other. Let $\varepsilon_{kl} \sim^{iid}$ Rademacher random variable $\varepsilon$ ($P[\varepsilon=\pm 1]=1/2$) which are also independent of A and B. The width is defined as (FACT 2.6)

$$W(\Gamma_X; \Phi_{A,B}):=E_H[sup\{\langle H,U\rangle: U \text{ in } \Gamma_X \text{ and } |U|_F=1] \quad (6.1)$$

where $$H:=m^{-1}\sum_{kl}\varepsilon_{kl}\Phi_{kl} =m^{-1}A^TEB \text{ and } E=[\varepsilon_{kl}] \quad (6.2)$$

and $\Gamma_X = D(|||.|||_1,X)$ in our applications.

Let $\boldsymbol{\varepsilon}_l$: $l=1,\ldots,m$ be column vectors of Rademacher matrix $[\varepsilon_{kl}]$, then $H=(\boldsymbol{h}_1,\ldots,\boldsymbol{h}_n)$ has its column vectors $\boldsymbol{h}_j=\sum_{l=1}^{m} m^{-1}B_{lj}A^T\boldsymbol{\varepsilon}_l \in R^n$. Notice that in problem $MP^{(F)}_{Y, A, B, \eta}$, $m^2$ is the measurement dimension.

### 6.1 Case 1: A and B are both Gaussian

**Lemma 6.1** Given $n$-by-$n$ matrix X=($\mathbf{x}_1,\dots,\mathbf{x}_n$) with $s$-sparse column vectors $\mathbf{x}_1,\dots,\mathbf{x}_n$, $r$ is cardinality of {$j$: $|\mathbf{x}_j|_1=max_k|\mathbf{x}_k|_1$}, $\Phi_{A,B}$, $\Gamma_X$, $W(\Gamma_X;\Phi_{A,B})$ are specified as the above and $A_{ki}\sim^{iid}B_{lj}\sim^{iid}N(0,1)$, then

$$W^2(\Gamma_X;\Phi_{A,B})\le 1+n^2-r(n-slog^2(cn^2r))$$

where $c$ is an absolute constant. Particularly, when $r = n$ ( i.e., X is sparse and $l_1$-column-flat ) then

$$W^2(\Gamma_X;\Phi_{A,B})\le 1+nslog^2(cn^3) \quad (6.3)$$

*Proof* We start with a similar inequality as that in FACT 2.4 (the proof is also similar) $W^2(\Gamma_X;\Phi_{A,B})\le E_H[inf\{|H-tV|_F^2: t>0, V \text{ in } \partial|||X|||_1\}]$. With the same specifications for V=($\lambda_1\xi_1,\dots,\lambda_n\xi_n$) as those in lemma 5.1, i.e.(w.l.o.g.) $\lambda_j\ge 0$ for $j$=1,…,$r$, $\lambda_1+\dots+\lambda_r=1$, $\lambda_j=0$ for $j\ge r+1$; $|\mathbf{x}_j|_1=max_k|\mathbf{x}_k|_1$ for $j$=1,…,$r$ and $|\mathbf{x}_j|_1 < max_k|\mathbf{x}_k|_1$ for $j\ge 1+r$; $\xi_j(i)=sgn(X_{ij})$ for $X_{ij}\ne 0$ and $|\xi_j(i)|\le 1$ for all $i$ and $j$. Let $\boldsymbol{h}_j\equiv\sum_{l=1}^m m^{-1}B_{lj}A^T\boldsymbol{\varepsilon}_l$, we have

$$W^2(\Gamma_X;\Phi_{A,B})\le E_{A,B,E}[inf_{t>0,\,\lambda j,\,\xi j \text{ specified as the above}}\sum_{j=1}^n|\sum_{l=1}^m m^{-1}B_{lj}A^T\boldsymbol{\varepsilon}_l-t\lambda_j\xi_j|_2^2]$$

$$=\sum_{j=r+1}^n E_{A,B,E}[|\boldsymbol{h}_j|_2^2]+E_{A,B,E}[inf_{t>0,\,\lambda j,\,\xi j \text{ specified as the above}}\sum_{j=1}^r|\boldsymbol{h}_j-t\lambda_j\xi_j|_2^2]$$

$$= \text{I} + \text{II}$$

The first and second terms are estimated respectively. The first term

$$\text{I}=\sum_{j=r+1}^n m^{-2}\sum_{l,k=1}^m E_B[B_{lj}B_{kj}]E_{A,E}[\boldsymbol{\varepsilon}_l^TAA^T\boldsymbol{\varepsilon}_k]=m^{-2}(n-r)\sum_{l,k=1}^m\delta_{lk}E_{A,E}[\boldsymbol{\varepsilon}_l^TAA^T\boldsymbol{\varepsilon}_l]=(n-r)n$$

To estimate II, for each $j$=1,…,$r$ let $S(j)$ be the support of $\mathbf{x}_j$ (so $|S(j)|\le$ s) and $\sim S(j)$ be its complimentary set, then

$$\sum_{j=1}^r|\boldsymbol{h}_j-t\lambda_j\xi_j|_2^2]=\sum_{j=1}^r|\boldsymbol{h}_{j|S(j)}-t\lambda_j\xi_{j|S(j)}|_2^2]+\sum_{j=1}^r|\boldsymbol{h}_{j|\sim S(j)}-t\lambda_j\xi_{j|\sim S(j)}|_2^2]$$

Notice that all components of $\xi_{j|S(j)}$ are $\pm 1$ and all components of $\xi_{j|\sim S(j)}$ can be any value in the interval [–1,+1]. Select $\lambda_1=\dots=\lambda_r=1/r$, let $\delta>0$ be arbitrarily small positive number and select $t=t(\delta)$ such that $P_{A,B,E}[|h|>t(\delta)/r]\le\delta$ where $h$ is a random scalar such that $h_j(i)\sim h$ and $i$ indicates the vector $\boldsymbol{h}_j$'s $i$-th component. For each $j$ and $i$ outside S($j$), set $\xi_j$'s component $\xi_j(i)=rh_j(i)/t(\varepsilon)$ if $|h_j(i)|\le t(\delta)/r$ and otherwise $\xi_j(i)=sgn(h_j(i))$, then $|\boldsymbol{h}_{j|\sim S(j)}-t\lambda_j\xi_{j|\sim S(j)}|_2^2=0$ when $|\boldsymbol{h}_{j|\sim S(j)}|_\infty< t(\delta)/r$ and notice the fact that for independent standard scalar Gaussian variables $a_l$, $b_l$ and Rademacher variables $\varepsilon_l$, $l$=1,…,$m$, there exists absolute constant $c$ such that for any $\eta>0$:

$$P[|m^{-1}\sum_{l,k=1}^m b_la_k\varepsilon_k|>\eta]<c\,exp(-\eta) \quad (6.4)$$

as a result, in the above expression $\delta$ can be $c\,exp(-t(\delta)/r)$ and:

$$E[|\boldsymbol{h}_{j|\sim S(j)}-t\lambda_j\xi_{j|\sim S(j)}|_2^2]=\int_0^\infty du P[|\boldsymbol{h}_{j|\sim S(j)}-t\lambda_j\xi_{j|\sim S(j)}|_2^2>u]=2\int_0^\infty duuP[|\boldsymbol{h}_{j|\sim S(j)}-t\lambda_j\xi_{j|\sim S(j)}|_2>u]$$

$$\le 2\int_0^\infty duuP[\text{There exists } (\boldsymbol{h}_{j|\sim S(j)}-t\lambda_j\xi_{j|\sim S(j)})\text{'s component with magnitude} > (n-s)^{-1/2}u]$$

$$\le 2(n-s)\int_0^\infty duuP[|h|-t(\delta)/r>(n-s)^{-1/2}u]$$

$$\le 2(n-s)\int_0^\infty duuexp(-((t(\delta)/r)+(n-s)^{-1/2}u))$$

$$\le C_0(n-s)^2exp(-(t(\delta)/r))\le C_0(n-s)^2\delta$$

where $C_0$ is an absolute constant. On the other hand $|\xi_{j|S(j)}|_2^2\le s$ for $j\ge 1+r$ so:

$$E_{A,B,E}[inf_{t>0,\,\lambda j,\,\xi j}\sum_{j=1}^r|\boldsymbol{h}_{j|S(j)}-t\lambda_j\xi_{j|S(j)}|_2^2]$$

$$\le E_{A,B,E}[\sum_{j=1}^r|\boldsymbol{h}_{j|S(j)}-t(\delta)\xi_{j|S(j)}/r|_2^2]$$

$$\le\sum_{j=1}^r E_{A,B,E}[m^{-2}|\sum_{l=1}^m B_{lj}(A^T\boldsymbol{\varepsilon}_l)_{|S(j)}|_2^2]+rst(\delta)^2/r^2$$

$$= rs(1+t(\delta)^2/r^2)$$

hence II $\le rs(1+t(\delta)^2/r^2)+nr\delta$. Combine all the above estimates we have:

$$W^2(\Gamma_X;\Phi_{A,B})\le \text{I}+\text{II}\le(n-r)n+rs(1+t(\delta)^2/r^2)+C_0n^2r\delta=n^2-r(n-s(1+t(\delta)^2/r^2)+C_0n^2r\delta$$

Substitute $t(\delta)/r$ with $log(c/\delta)$ we get, for any $\delta>0$:

$$W^2(\Gamma_X;\Phi_{A,B})\le n^2-r(n-s(1+log^2(c/\delta))+C_0n^2r\delta$$

In particular, let $\delta=1/C_0n^2r$ then $W^2(\Gamma_X;\Phi_{A,B})\le n^2-r(n-s(1+log^2(cn^2r))+1$. □

To apply FACT 2.6 completely, now estimate the lower bound of $Q_\xi(\Gamma_X;\Phi_{A,B}):=inf\{P[|\langle\Phi_{kl},U\rangle|\ge\xi]$:

U in $\Gamma_X$ and $|U|_F=1\}$ for some $\xi > 0$ where $<\Phi_{kl},U>=\sum_{ij}A_{ki}U_{ij}B_{lj}$ for each $1\le k,l\le m$, $A_{ki} \sim^{iid} B_{lj} \sim^{iid} N(0,1)$. The estimate is independent of the indices $k$ and $l$, so we give a general and notational simplified statement on it.

**Lemma 6.2** Let $M_{ij} = a_i b_j$ where $a_i \sim^{iid} b_j \sim^{iid} N(0,1)$ and all random variables are independent each other, there exists an positive absolute constant $c$ such that $inf\ \{P[|<M,U>| \ge 1/\sqrt{2}]$: U in $\Gamma_X$ and $|U|_F=1\} \ge c$, as a result $Q_{1/\sqrt{2}}(\Gamma_X; \Phi_{A,B}) \ge c$.

*Proof* By the second moment inequality $P[Z \ge \xi] \ge (E[Z] - \xi)_+^{\ 2}/E[Z^2]$ for any non-negative r.v. Z and any $\xi > 0$. Set $Z = |<M,U>|^2$ and $\xi = E[|<M,U>|^2]/2$, we get:

$$P[|<M,U>|^2 \ge E[|<M,U>|^2]/2] \ge E[|<M,U>|^2]^2/4E[|<M,U>|^4] \qquad (6.5)$$

To estimate the upper bound of $E[|<M,U>|^2]$, let $U=\sum_j \lambda_j \mathbf{u}_j \mathbf{v}_j$ be U's singular value decomposition, $\mathbf{u}_i^T\mathbf{u}_j=\mathbf{v}_i^T\mathbf{v}_j=\delta_{ij}$, $\lambda_j>\mathbf{0}$ for each $j$. Notice that $M=\mathbf{ab}^T$ where $\mathbf{a}\sim\mathbf{b}\sim N(0, I_n)$ and independent each other, then $<M,U> = \mathbf{a}^T U\mathbf{b} = \sum_j \lambda_j \mathbf{a}^T\mathbf{u}_j\mathbf{v}_j^T\mathbf{b}$ where $\mathbf{a}^T\mathbf{u}_i \sim \mathbf{v}_j^T\mathbf{b} \sim N(0,1)$ and independent each other, hence $E[|<M,U>|^2] = \sum_j \lambda_j^2 E[|\mathbf{a}^T\mathbf{u}_j|^2]E|\mathbf{v}_j^T\mathbf{b}|^2] = \sum_j \lambda_j^2 = |U|_F^2=1$ for U in the assumption.

On the other hand by Gaussian hypercontractivity we have

$$(E[|<M,U>|^4])^{1/4} \le C_0(E[|<M,U>|^2])^{1/2} = C_0$$

In conclusion $P[|<M,U>|^2 \ge 1/2] = P[|<M,U>|^2 \ge E[|<M,U>|^2]/2] \ge c$ for U: $|U|_F^2=1$. □

Combing lemma 6.1, 6.2 and FACT 2.6, we obtain the general result in the following.

**Theorem 6.1** Suppose $A_{ki} \sim^{iid} B_{lj} \sim^{iid} N(0,1)$ and independent each other, $X\in \sum_s^{n\times n}$ is a columnwise $s$-sparse and $l_1$-column-flat signal, $\mathbf{Y} = AXB^T + E \in R^{m\times m}$ with measurement errors bounded by $|E|_F^2 \le \eta$, $X^*$ is the minimizer of the problem $MP^{(F)}_{Y, A,B, \eta}$. If

$$m \ge t + 4\sqrt{2}\eta/\delta + C_1(ns)^{1/2}log(C_2 n^3)$$

where $C_i$'s are absolute constants, then $P[|X^*-X|_F \le \delta] \ge 1 - exp(-t^2/2)$, i.e., such X can be reconstructed robustly with respect to the error norm $|X^*-X|_F$ with high probability by solving $MP^{(F)}_{Y,A,B,\eta}$. □

## 6.2 Case 2: A and B are both sub-Gaussian

**Lemma 6.3** Given $n$-by-$n$ matrix $X=(\mathbf{x}_1,\ldots,\mathbf{x}_n)$ with $s$-sparse column vectors $\mathbf{x}_1,\ldots,\mathbf{x}_n$, $r$ is cardinality of $\{j: |\mathbf{x}_j|_1=max_k|\mathbf{x}_k|_1\}$, $\Phi_{A,B}$, $\Gamma_X$, $W(\Gamma_X; \Phi_{A,B})$ are specified as before, $A_{ki} \sim^{iid}$ *Sub-Gaussian* distribution and $B_{lj} \sim^{iid}$ *Sub-Gaussian* distribution with $\psi_2$-norms $\sigma_A$, $\sigma_B$ respectively, then

$$W^2(\Gamma_X; \Phi_{A,B}) \le \sigma_A^2\sigma_B^2(1+n^2-r(n-slog^2(Cn^2r))) \qquad (6.6)$$

where $C$ is an absolute constant. Particularly, when $r = n$ then

$$W^2(\Gamma_X; \Phi_{A,B}) \le \sigma_A^2\sigma_B^2(1 + nslog^2(Cn^3)) \qquad (6.7)$$

The proof of this lemma is logically the same as the proof of lemma 6.1, the only difference is about the distribution tail of the components of vectors $\boldsymbol{h}_j \equiv \sum_{l=1}^m m^{-1}B_{lj}A^T\boldsymbol{\varepsilon}_l$ which $\sim^{iid}$ $h \equiv m^{-1}\sum_{l,k=1}^m b_l a_k \varepsilon_k$ with independent scalar sub-Gaussian variables $a_l$, $b_l$ and Rademacher variables $\varepsilon_l$, $l=1,\ldots,m$. This auxiliary result is presented in the following lemma:

**Lemma 6.4** For independent scalar zero-mean *sub-Gaussian* variables $a_l$, $b_l$ and Rademacher variables $\varepsilon_l$, $l=1,\ldots,m$, let $\sigma_A\equiv max_l|a_l|_{\psi 2}$, $\sigma_B\equiv max_l|b_l|_{\psi 2}$ ($|.|_{\psi 2}$ denotes a *sub-Gaussian* variable's $\psi_2$-norm), then there exists absolute constant $c$ such that for any $\eta > 0$:

$$P[|\ h\ | > \eta] < 2\ exp(-c\eta/\sigma_A\sigma_B) \qquad (6.8)$$

*Proof* Notice that $a_k\varepsilon_k$ is zero-mean *sub-Gaussian* variable with $|a_k\varepsilon_k|_{\psi 2}=|a_k|_{\psi 2}$, for $b=m^{-1/2}\sum_{1\le l\le m}b_l$ and $a=m^{-1/2}\sum_{1\le k\le m}a_k\varepsilon_k$ we have $|b|_{\psi 2} \le Cm^{-1/2}(\sum_l|b_l|_{\psi 2}^2)^{1/2} \le C\sigma_B$ and $|a|_{\psi 2} \le Cm^{-1/2}(\sum_l|a_k|_{\psi 2}^2)^{1/2} \le C\sigma_A$ where $C$ is an absolute constant. Furthermore, because the product of two *sub-Gaussian* variables $a$ and $b$ is

*sub-Exponential* and its $\psi_1$-norm $|ba|_{\psi 1} \leq |b|_{\psi 2}\,|a|_{\psi 2} \leq C^2\sigma_A\sigma_B$, $h \equiv m^{-1}\sum_{l,k=1}^{m} b_l a_k \varepsilon_k = ab$ has its distribution tail $P[|\,h\,| > \eta] < 2exp(-c\eta/\sigma_A\sigma_B)$ where $c$ is an absolute constant.

*Proof of lemma 6.3* With the same logic as in the proof of lemma 6.1 and based-upon lemma 6.4, the auxiliary parameter $\delta$ in the argument can be $2exp(-ct(\delta)/r\sigma_A\sigma_B)$ and equivalently $t(\delta)/r = \sigma_A\sigma_B log(2/\delta)$ which derives the final result. □

**Lemma 6.5** Let $M_{ij} = a_i b_j$ where $a_i \sim^{iid}$ *Sub-Gaussian* distribution, $b_j \sim^{iid}$ *Sub-Gaussian* distribution and all are independent each other, then $Q_{1/\sqrt{2}}(\Gamma_X;\ \Phi_{A,B}) \geq c$ where $c$ is a constant only dependent on $\sigma_A\sigma_B$. □

Lemma 6.5's proof is almost the same as that of lemma 6.2. Now the following result can be obtained directly by combing lemmas 6.3-6.5 and FACT 2.6:

**Theorem 6.2** Suppose random matrices A, B are independent each other, $A_{ki} \sim^{iid}$ *Sub-Gaussian* distribution, $B_{lj} \sim^{iid}$ *Sub-Gaussian* distribution, each with $\psi_2$-norm $\sigma_A$ and $\sigma_B$. Let $X \in \sum_s^{n\times n}$ be a columnwise $s$-sparse and $l_1$-column-flat signal, $\mathbf{Y} = AXB^T + E \in R^{m\times m}$ with $|E|_F^2 \leq \eta$, $X^*$ be the minimizer of $MP^{(F)}_{Y,A,B,\eta}$. If

$$m \geq t + 4\sqrt{2}\eta/\delta + C_1\sigma_A\sigma_B(ns)^{1/2}log(C_2n^3)$$

where $C_i$'s are absolute constants, then $P[|X^*-X|_F \leq \delta] \geq 1 - exp(-t^2/2)$, i.e., such X can be reconstructed robustly with respect to the error norm $|X^*-X|_F$ with high probability by solving $MP^{(F)}_{Y,A,B,\eta}$. □

## 7 Conclusions, Some Extensions and Future Works

In this paper we investigated the problem of reconstructing $n$-by-$n$ column-wise sparse and $l_1$-column-flat matrix signal $X=(\mathbf{x}_1,\ldots,\mathbf{x}_n)$ via convex programming with the regularizer $|||X|||_1 := max_j|\mathbf{x}_j|_1$ where $|.|_1$ is the $l_1$-norm in vector space. In the first part (sec.3 and sec.4), the most important conclusions are about the general conditions to guarantee uniqueness, value-robustness, support stability and sign stability in signal reconstruction. In the second part (sec.5 and sec.6) we took the convex geometric approach in random measurement setting and established bounds on dimensions of measurement spaces for robust reconstruction in noise. For example, typical results show that the signal complexity (width) encoded by the regularizer $|||X|||_1$ is determined by two structural parameters, i.e., the maximum number $s$ of nonzero entries in each column (column-wise sparsity )and the number $r$ of columns which $l_1$-norms are maximum among all columns ($l_1$-column-flatness). The signal's complexity decreases with small $s$ and large $r$. In particular, when $r = n$ (i.e., sparse and flat signal) the condition reduces to $m \geq t + 4\sqrt{2}\eta/\delta + C_1\sigma_A\sigma_B(ns)^{1/2}log(C_2n^3)$.

Based on the methods and results obtained in this paper, we can make some straightforward extensions to other similar problems. The first extension is about reconstructing row-wise sparse and $l_1$-row-flat matrix signals via convex programming with regularizer $|||X^T|||_1 := max_i|\overline{\mathbf{x}_i}|_1$ where $\overline{\mathbf{x}_i}$ is matrix signal X's $i$-th row. In this case all the obtained estimates and conclusions remain the same, e.g., the signal's complexity is determined by the maximum number $s$ of nonzero entries in each row and the number $r$ of rows which $l_1$-norms are maximum among all rows. The signal's complexity decreases with small $s$ and large $r$.

The second extension is for reconstructing the matrix signal which has both row-wise and column-wise sparsity and $l_1$-flatness. We can use $F(X) := max(|||X|||_1, |||X^T|||_1)$ as the regularizer. Note that when $|||X|||_1 \geq |||X^T|||_1$ then $F(X) = |||X|||_1$ so in this case $W(\Gamma_X;\ \Phi)$ is the width determined by regularizer $|||X|||_1$, otherwise $W(\Gamma_X;\ \Phi)$ is the width determined by regularizer $|||X^T|||_1$, both have the upper bound of $\overline{W}(r_1,s_1,n)$ and $\overline{W}(r_2,s_2,n)$ in the form as that in lemma 5.1 or lemma 6.1 where $(r_1,s_1)$ and $(r_2,s_2)$ are the matrix signal's row and column structural parameters. As a result, the sufficient condition for robust reconstruction is $m \geq max(\overline{W}^2(r_1,s_1,n),\ \overline{W}^2(r_2,s_2,n))$ via $MP_{y,\Phi,\eta}$ or $m^2 \geq max(\overline{W}^2(r_1,s_1,n),\ \overline{W}^2(r_2,s_2,n))$ via $MP_{Y,A,B,\eta}$.

The third extension is for reconstructing the matrix signal with the general linear measurement Y = $\Phi_{A,B}$(X) + E = $\sum_{\mu=1}^{L} A_\mu X B_\mu^T$ + E where $|E|_F \leqq \eta$. In this case (see (6.2)) $H = \sum_{\mu=1}^{L} H_\mu$ where $H_\mu := m^{-1} A_\mu^T E_\mu B_\mu$ and $E_\mu = [\varepsilon_{kl}^{(\mu)}]$ is the matrix with independent Rademacher entries. Suppose matrices $A_\mu$'s and $B_\mu$'s are both Sub-Gaussian, in this case (see (4.1)) the width

$$W(\Gamma_X; \Phi_{A,B}) = E_H[sup\{\sum_{\mu=1}^{L} <H_\mu, U>: U \text{ in } \Gamma_X \text{ and } |U|_F=1]$$

$$\leq C_0 L^{1/2} max_\mu E_H[sup\{<H_\mu, U>: U \text{ in } \Gamma_X \text{ and } |U|_F=1]$$

and each $E_H[sup\{|<H_\mu,U>|^2: U$ in $\Gamma_X$ and $|U|_F=1]$ has the upper bound estimate $\sigma_A^2\sigma_B^2(n^2-r(n-slog^2(C_1n^2r))$ according to lemma 4.3, finally we get the width estimate $W^2(\Gamma_X; \Phi_{A,B}) \leq C_0 L\sigma_A^2\sigma_B^2(n^2-r(n-slog^2(C_1n^2r))$. As a result, a sufficient condition for column-wise sparse and $l_1$-column-flat signal X to be reconstructed via solving $MP^{(F)}_{Y,A,B,\eta}$ robustly with respect to *Frobenius* norm $|.|_F$ from linear measurement Y = $\sum_{\mu=1}^{L} A_\mu X B_\mu^T$+E where $|E|_F \leqq \eta$ is $m \geq t + 4\sqrt{2}\eta/\delta + C_1 L\sigma_A\sigma_B(ns)^{1/2}log(C_2n^3)$.

This paper is only focused on basic theoretical analysis. In subsequent papers, the algorithms to solve the $|||.|||_1$-optimization problems (e.g., generalized inverse scale space algorithms, etc.), related numeric investigations and applications (e.g.,in radar space-time waveform analysis)will be further investigated[14-15].